\documentclass{article}
\usepackage[a4paper,top=4cm,bottom=4cm,left=2.5cm,right=2.5cm,bindingoffset=5mm]{geometry}
\usepackage{xcolor}

\usepackage[font=small,format=default,indention=2em,labelsep = period]{caption} 

\usepackage{enumitem}
\usepackage{subfig}
\usepackage{amssymb,bm,amsmath}
\usepackage{graphicx,algorithm,algorithmic,mathtools}	
			
\newcommand{\R}{\mathbb R}
\newcommand{\dt}{\Delta t}
\newcommand{\nt}{n_{\textup{f}}}
\newcommand\norm[1]{\left\lVert#1\right\rVert}
\newcommand{\init}{\textup{i}}
\newcommand{\fin}{\textup{f}}
\newcommand{\pdata}{{\bf P^{\mathrm{data}}}}
\newcommand{\pdmd}{{\bf {\widetilde{P}}}^{DMD}}
\DeclarePairedDelimiter{\floor}{\lfloor}{\rfloor}
\DeclareRobustCommand{\bbone}{\text{\usefont{U}{bbold}{m}{n}1}}  
\DeclareMathOperator{\dimvett}{dim}

\usepackage{titlesec}
\titleformat*{\section}{\large\bfseries}
\titleformat*{\subsection}{\normalsize\bfseries}
\let\oldparagraph=\paragraph
\renewcommand\paragraph[1]{\oldparagraph{#1.}}

\numberwithin{equation}{section}
\usepackage{ntheorem}   
\theoremstyle{plain}
\theoremheaderfont{\bfseries}
\theorembodyfont{\normalsize}
\theoremseparator{.}
\newtheorem{remark}{Remark}
\numberwithin{remark}{section}

\title{\large\textbf{Understanding Mass Transfer Directions via Data-Driven Models \\with Application to Mobile Phone Data}}
\author{\normalsize{Alessandro Alla}\thanks{Departamento de Matem\'atica, PUC-Rio, Rio de Janeiro, BR (\href{mailto:alla@mat.puc-rio.br}{alla@mat.puc-rio.br}).}
\and \normalsize{Caterina Balzotti}\thanks{Dipartimento di Scienze di Base e Applicate per l'Ingegneria, Sapienza Universit\`a di Roma, Rome, IT (\href{mailto:caterina.balzotti@sbai.uniroma1.it}{caterina.balzotti@sbai.uniroma1.it}).}
\and \normalsize{Maya Briani}\thanks{Istituto per le Applicazioni del Calcolo ``M.\ Picone'', Consiglio Nazionale delle Ricerche, Rome, IT (\href{mailto:m.briani@iac.cnr.it}{m.briani@iac.cnr.it}, \href{mailto:e.cristiani@iac.cnr.it}{e.cristiani@iac.cnr.it}).}
\and \normalsize{Emiliano Cristiani}\footnotemark[3]}
\date{\vspace{-0.5cm}}

\begin{document}

\maketitle

\begin{abstract}
The aim of this paper is to solve an inverse problem which regards a mass moving in a bounded domain. We assume that the mass moves following an unknown velocity field and that the evolution of the mass density can be described by a partial differential equation, which is also unknown. The input data of the problems are given by some snapshots of the mass distribution at certain times, while the sought output is the velocity field that drives the mass along its displacement.
To this aim, we put in place an algorithm based on the combination of two methods: first, we use the Dynamic Mode Decomposition to create a mathematical model describing the mass transfer; second, we use the notion of Wasserstein distance (also known as earth mover's distance) to reconstruct the underlying velocity field that is responsible for the displacement.
Finally, we consider a real-life application: the algorithm is employed to study the travel flows of people in large
populated areas using, as input data, density profiles (i.e.\ the spatial distribution) of people in given areas at different time instants. This kind of data are provided by the Italian telecommunication company TIM and are derived by mobile phone usage.
\end{abstract}

\begin{description}
\item[\textbf{Keywords.}] Data-driven methods, dynamic mode decomposition, Wasserstein distance, earth mover's distance, cellular data, presence data.
\item[\textbf{Mathematics Subject Classification.}]  37C10, 35R30, 76D55.
\end{description}

\section{Introduction}

In this paper we aim at solving an inverse problem which regards a mass moving in a bounded domain with finite velocity. We assume that the mass moves following an unknown velocity field and that the evolution of the mass density can be described by an unknown PDE. The input data of the problem are given by some snapshots of the mass distribution at certain times, while the sought information is the velocity field that drives the mass along its displacement.

The basic idea can be summarized as follows: given two snapshots of the mass distribution at two instants of time, we want to understand where each portion of the mass (which is assumed to be conserved from one instant of time to the other) is transported from/to, i.e.\ how the first spatial distribution is rearranged in the second one. 
The goal is pursued by computing a numerical approximation of the Wasserstein distance (also known as earth mover's distance or Mallows distance) between the two consecutive density profiles, specifying a suitable cost function which measures the ``energy'' consumed by the system for moving forward. The computation of the Wasserstein distance gives, as a by-product, the minimum-cost mass flow from the first to the second configuration, i.e.\ how the mass distributes in space and time.

Despite the good preliminary results obtained by the above described Wasser\-stein-based approach \cite{balzotti2018IFAC, zhu2018IJGIS}, the algorithm is found to be excessively expensive both in terms of CPU time and memory requirements. This fact strongly restricts the applicability of the method. 
To fix this, in this paper we propose to couple the method with the Dynamic Mode Decomposition (DMD): a data-driven technique that takes in input the snapshots of the mass distribution and returns an analytical approximation of the dynamics underlying the mass transfer. 
More precisely, it provides a system of ODEs which describes the evolution of the mass in any point of the domain. Solving the ODEs, we are able to recover the mass distribution \emph{at any time}, thus increasing at will the number of available snapshots or, analogously, decreasing at will the time frame between them. 
Controlling the time frame between two consecutive snapshots is the key to simplify the computation of the Wasserstein distance and makes the computation of the flows feasible even on large domains. 

Finally, a real-world application of the proposed methodology is illustrated. We are interested in inferring activity-based human mobility flows from mobile phone data. 
We assume that mobile devices are not singularly tracked, but their logs are aggregated in order to obtain the total number of users in a given area. In this way we get the density profiles (i.e.\ the spatial distribution) of people in a given area at various instants of time. The dataset we have at our disposal is provided by the Italian telecommunication company TIM. The time frame between two consecutive snapshots is 15 minutes.

As before, the goal is to ``assign a direction'' to the presence data. In fact, the mere representation of time varying density of people clearly differentiate attractive from repulsive or neutral areas but does not provide any information about the directions of flows of people. 
In other words, we are interested in a ``where-from-where-to'' type of information, which reveals travel flows and patterns of people providing a sort of \emph{origin-destination} matrix.

\paragraph{Paper organization}  The paper is organized as follows. 
In Section \ref{sec:background} we present the DMD method and the Wasserstein distance. In Section \ref{sec:coupling} we show how to couple those methods to obtain an efficient algorithm.
In Section \ref{sec:TIM} we apply the proposed approach to real-life data. In  Section \ref{sec:conclusions} we draw our conclusions. Finally, in Appendix \ref{sec:toyexampleDMD} we apply the DMD method to the viscous Burgers' equation.

\section{Mathematical background}\label{sec:background}
In this section we recall the building blocks of the methodology proposed in the paper, namely the DMD method used to build a data-driven model, and the Wasserstein distance to determine the transport map driving the moving mass.

\subsection{DMD method}\label{sec:introDMD}

DMD is a data-driven method capable of providing accurate assessments of the spatio-temporal coherent structures in a given complex system, or short-time future estimates of such a system. Although a complete list of references for DMD goes beyond the scopes of this work, we would like to mention \cite{BN15, HH18b, HH18a, S10,TCBN14}. In the current work we use the DMD algorithm as in \cite{TRLBK14}. 

To begin, we suppose to have a set of data $\mathcal{X}=\{{\bf y}(t_0),\ldots,{\bf y}(t_{\nt})\}$ for some time instances $\{t_j\}_{j=0}^{\nt}$ with ${\bf y}(t_j)\in\R^N, j=0, \ldots, \nt$ and $\Delta t = t_{j+1}-t_j$ for $j=0,\dots,\nt-1$. The goal of the method is to build a mathematical model upon the dataset $\mathcal{X}\in\R^{N\times (\nt+1)}$.
The DMD procedure thus constructs the approximate linear evolution $\widehat{\bf y}(t)$ for the dataset $\mathcal{X}$ exploiting its low-rank structure:
\begin{equation}
\frac{d \widehat{\bf y}}{dt} = \widehat{\bf A} \widehat{\bf y} \label{eq:dA}
\end{equation}
where $\widehat{\bf A}\in\R^{N\times N}$ is unknown, $\widehat{\bf y}(0) = \widehat{\bf y}_0$, and the solution has the form
\begin{equation}
\widehat{\bf y}(t)=  \sum_{i=1}^r \beta_i {\boldsymbol{\psi}}_i \exp(\omega_i t) \, ,
\label{eq:omegaj}
\end{equation}
where $r<N, {\boldsymbol{\psi}}_i$ and $\omega_i$ are the eigenvectors  and eigenvalues of the unknown matrix $\widehat{\bf A}$. 
 The coefficients $\beta_i$ of the vector ${\boldsymbol \beta}$ can be determined from the initial data.  For example, at $t=t_0$ we have ${\bf y}(t_0)={\bf y}_0$ so that \eqref{eq:omegaj} gives ${\boldsymbol \beta}={\bf \Psi}^\dag {\bf y}_0$, where ${\bf \Psi}$ is a matrix comprised of the DMD modes $\boldsymbol{\psi}_i$ and $\dag$ denotes the Moore-Penrose pseudo-inverse.
To compute the matrix $\widehat{\bf A}$, we first split the dataset into two snapshot matrices
%
\begin{equation}\label{inp_out}
  {\bf Y} \!=\! \begin{bmatrix}
\vline & \vline & & \vline \\
{\bf y}(t_0) & {\bf y}(t_1) & \cdots & {\bf y}(t_{\nt-1})\\
\vline & \vline & & \vline
\end{bmatrix}, \hspace{0.1in} {\bf Y}' \!=\! \begin{bmatrix}
\vline & \vline & & \vline \\
{\bf y}(t_1) & {\bf y}(t_2) & \cdots & {\bf y}(t_{\nt})\\
\vline & \vline & & \vline
\end{bmatrix}
\end{equation}
%
and suppose the following linear relation hold true:
\begin{equation}
{\bf Y}'={\bf AY},
\end{equation}
where ${\bf A}:=\exp{(\widehat{\bf A}\Delta t}$).

Specifically, we assume that ${\bf y}(t_j)$ is an initial condition to obtain ${\bf y}(t_{j+1})$, i.e. its corresponding output after some prescribed evolution time $\Delta t>0$. Thus, the DMD method computes the best linear operator ${\bf A}$ relating to the matrices above:
\begin{equation}
  {\bf A} = {\bf Y}' {\bf Y}^\dag.
  \label{eq:newDMD}
\end{equation}
 We refer to ${\bf Y}$ and ${\bf Y}'$ as input and output snapshot matrices respectively.

The DMD algorithm aims at optimally constructing the matrix ${\bf A}$ so
that the error between the true and approximate solution is small in a least-square sense, i.e.  
$\| {\bf y}(t) - \widehat{\bf y}(t) \| \ll 1$. 
Of course, the optimality of the approximation holds only over the sampling window where ${\bf A}$ is constructed, but the approximate solution
can be used to make future state predictions, and to decompose the dynamics into various time-scales.

The matrix ${\bf A}$ is often highly ill-conditioned and when the state dimension $n$ is large, the aforementioned matrix may be even intractable to analyze directly. Instead, DMD circumvents the eigen-decomposition of ${\bf A}$ by considering a low rank representation in terms of a matrix $\widehat{\bf A}$ projected with the Proper Orthogonal Decomposition (POD). Although the description of the POD method goes beyond the scopes of this paper, we recall that the POD projection solves the following optimization problem
\begin{equation}\label{pbmin}
\min_{ {\boldsymbol{\varphi}}_1,\ldots,{\boldsymbol{\varphi}}_r\in\R^n} \sum_{j=0}^{\nt-1} \left\|{\bf y}(t_j)-\sum_{i=1}^r \langle {\bf y}(t_j),{\boldsymbol{\varphi}}_i\rangle{\boldsymbol{\varphi}}_i\right\|^2\quad \mbox{such that }\langle {\boldsymbol{\varphi}}_i,{\boldsymbol{\varphi}}_j\rangle=\delta_{ij},
\end{equation}
where $\{{\boldsymbol\varphi}_i\}_{i=1}^r$ are the POD projectors.
The solution of the optimization problem \eqref{pbmin} is obtained by means of a Singular Value Decomposition (SVD) of the dataset ${\bf Y}$, where the first $r\ll N$ columns of the left singular eigenvectors are the required POD basis.  We refer to \cite{Vol11} for a complete description of the POD method. We also mention that POD method is equivalent to Principal Component Analysis (PCA) or Karhunen-Lo\'eve expansion in other contexts (see e.g. \cite{J02,P01}).  

The exact DMD algorithm proceeds as follows \cite{TRLBK14}:
first, we collect data ${\bf Y, Y'}$ as in \eqref{inp_out} and compute the reduced, or economy, singular value decomposition of {\bf Y}
$${\bf Y}={\bf U}{\bf \Sigma}{\bf V}^T.$$ 

We note that the use of the economy SVD is suggested since the matrices ${\bf Y}, {\bf Y}'\in\R^{n\times \nt}$ with $N\gg \nt$ and that the economy ${\bf U}\in\R^{N\times r}$ is sufficient to provide the same approximation of the regular SVD given the limited amount of snapshots. \\ 
Furthermore, the DMD typically makes use of low-rank structure so that the total number of modes, $r\ll N$, allows for dimensionality reduction of the
dynamical system.
Then, we compute the least-squares fit ${\bf A}$ that satisfies ${\bf Y}'={\bf AY}$ and project onto the POD modes ${\bf U}$.
Specifically, the Moore-Penrose pseudo-inverse of ${\bf Y}$ allows us to compute ${\bf A}={\bf Y}'{\bf Y}^\dag$, where the Moore-Penrose algorithm provides the least-square fitting procedure. In terms of its low-rank projection, this yields
$${\bf \widehat{A}}={\bf U}^T {\bf AU}={\bf U}^T{\bf Y}'{\bf V}{\bf \Sigma}^{-1},$$
and then, we compute the eigen-decomposition of ${\bf \widehat{A}}\in\R^{r\times r}$
$${\bf\widehat A}{\bf W}={\bf W \Lambda},$$
where ${\bf \Lambda}$ are the DMD eigenvalues.
Finally, the DMD modes ${\bf \Psi}^{\mbox{\tiny DMD}}$ are given by
\begin{equation}\label{dmd_basis}
{\bf\Psi}^{\mbox{\tiny DMD}}={\bf Y}'{\bf V}{\bf \Sigma}^{-1}{\bf W}.
\end{equation}
The algorithm is summarized in Algorithm \ref{Alg_DMD}.
\begin{algorithm}[H]
\caption{Exact DMD}
\label{Alg_DMD}
\begin{algorithmic}[1]
\REQUIRE Snapshots $\{{\bf y}(t_0),\ldots,{\bf y}(t_{\nt})\}$, Time step $\Delta t$.
\STATE Set ${\bf Y}=[{\bf y}(t_0),\dots, {\bf y}(t_{\nt-1})]$ and ${\bf Y'}=[{\bf y}(t_1),\dots, {\bf y}(t_{\nt})]$.
\STATE Compute the reduced SVD of rank $r$ of ${\bf Y}$, ${\bf Y}={\bf U}{\bf \Sigma}{\bf  V}^T$.
\STATE Define $\widehat{{\bf A}}:={\bf U}^T{\bf Y}'{\bf V}{\bf \Sigma}^{-1}$.
\STATE Compute eigenvalues and eigenvectors of $\widehat{{\bf A}} {\bf W}={\bf W}{\bf \Lambda}$.
\STATE Set ${\bf \Psi}^{\mbox{\tiny DMD}}={\bf \Lambda}^{-1}{\bf Y}'{\bf V}{\bf \Sigma}^{-1}{\bf W}$. 
\STATE Set $\omega_i = \frac{\log(\lambda_i)}{\Delta t}$ for \eqref{eq:omegaj}
\end{algorithmic}
\end{algorithm}
In Appendix \ref{sec:toyexampleDMD} we provide a numerical experiment to show the effectiveness of the DMD method. We compute the data from a nonlinear PDE, i.e. the viscous Burgers' equation.

\subsection{Wasserstein distance and optimal mass transfer problem}\label{sec:wassersteindistance}
The notion of Wasserstein distance is strictly related to the Monge--Kanto\-rovich optimal mass transfer problem \cite{villani2008SSBM}, which can be easily explained as follows: given a sandpile with mass distribution $\rho^\init$ and a pit with equal volume and mass distribution $\rho^\fin$, find a way to minimize the cost of transporting sand into the pit. The cost for moving mass depends on both the distance from the point of origin to the point of arrival and the amount of mass is moved along that path. We are interested in minimizing this cost by finding the optimal path to transport the mass from the initial to the final configuration.

Given two density functions $\rho^\init, \rho^\fin:\Omega\to\R$ for some bounded $\Omega\subset\R^n$ such that $\int_{\R^{n}}\rho^{\init}=\int_{\R^{n}}\rho^{\fin}$, we define the $L^p$-Wasserstein distance between $\rho^\init$ and $\rho^\fin$ as
\begin{equation}\label{def:WassDist}
W_p(\rho^\init,\rho^\fin)=\bigg(\min_{T\in\mathcal{T}}\int_{\Omega}c(\bm{\xi},T(\bm{\xi}))^p \, \rho^\init(\bm{\xi})d\bm{\xi}\bigg)^{\frac{1}{p}}
\end{equation}
where 
\begin{equation*}
\mathcal{T}:=\Biggr\{T\colon\Omega\to\Omega \, : \, \int\displaylimits_B  \rho^\fin(\bm{\xi})d\bm{\xi} =\int\displaylimits_{\{\bm{\xi}:T(\bm{\xi})\in B\}} \rho^\init(\bm{\xi})d\bm{\xi}, \quad
\forall \, B\subset\Omega\Biggr\}
\end{equation*}
and $c:\Omega\times\Omega\to\R$ is a given cost function, which defines the cost of transferring a unit mass between any two points in $\Omega$.
Note that $\mathcal{T}$ is the set of all possible maps which transfer the mass from one configuration to the other. 

It is important to note here that we are not really interested in the actual value of the Wasserstein distance $W_p$, instead we look for the \emph{optimal map} $T^*$ which realizes the $\arg\min$ in \eqref{def:WassDist}, and represents the paths along which the mass is transferred.

\subsection{Numerical approximation of the Wasserstein distance}\label{sec:hitchcock}
A direct numerical approximation of definition \eqref{def:WassDist} is unfeasible, but a discrete approach is still possible. Indeed, we can resort to classical problems (see Hitchcock's paper \cite{hitchcock1941SAM}) and
methods (see e.g., \cite[Sec.\ 6.4.1]{santambrogio2015SPRINGER} and \cite[Chap.\ 19]{sinha2005ELSEVIER}) to recast the original mass transfer problem in the framework of linear programming (LP). We also refer to \cite{briani2017CMS} for a recent application of this methodology to traffic flow problems. 

The idea consists in approximating the set $\Omega$ with a structured grid with $N$ square cells $C_1,\ldots,C_N$, as it is commonly done for the numerical approximation of PDEs. We denote by $\Delta x$ the length of each side of the cells. Then, we define a graph $\mathcal{G}$ whose nodes coincide with the centers of the $N$ cells.  
Graph's edges are defined in such a way that each node is directly connected with each other, including itself.

Introducing a numerical error (controlled by the parameter $\Delta x$), we are allowed to assume that $\forall j=1,\ldots, N$ all the mass distributed in the cell $C_j$ is concentrated in its center, i.e.\ in a node of the graph. We come up with an initial mass $m^\init_j:=\int_{C_j}\rho^\init dx$ and a final mass $m^\fin_j:=\int_{C_j}\rho^\fin dx$, for $j=1,\dots,N$, distributed on the graph nodes.
Now, we simply aim at optimally rearranging the first mass into the second one moving it among the graph's nodes. 

We denote by  $c_{jk}$ the cost to transfer a unit mass from node $j$ to node $k$, and by $x_{jk}$ the (unknown) mass moving from node $j$ to node $k$. The problem is then formulated as
\begin{equation*} 
\text{minimize }\mathcal{H}:= \sum_{j,k=1}^N c_{jk}x_{jk}
\end{equation*}
subject to
\[\displaystyle\sum_k x_{jk}=m_j^\init \quad \forall j, \quad
\displaystyle\sum_j x_{jk}=m_k^\fin \quad \forall k \quad\text{and}\quad x_{jk}\geq 0.\]
Defining
\begin{flalign*}
&{\bf x} = (x_{11}, x_{12},\dots,x_{1N},x_{21},\dots,x_{2N},\dots,x_{N1},\dots,x_{NN})^\textup{T},\\
&{\bf c}= (c_{11}, c_{12},\dots,c_{1N},c_{21},\dots,c_{2N},\dots,c_{N1},\dots,c_{NN})^\textup{T},\\
&{\bf b} = (m^\init_1,\dots,m^\init_N,m^\fin_1,\dots,m^\fin_N)^\textup{T},
\end{flalign*}
and the $2N\times N^2$ matrix
\[\bf M =\begin{bmatrix}
\boldsymbol\bbone_N & 0 & 0 & \dots & 0\\[0.3em]
0 & \boldsymbol\bbone_N & 0 & \dots & 0\\[0.3em]
0 & 0 & \boldsymbol\bbone_N &  \dots & 0\\[0.3em]
\vdots & \vdots & \vdots & \ddots & \vdots\\[0.3em]
0 & 0 & 0 & \dots & \boldsymbol\bbone_N\\[0.3em]
{\bf I}_N &{\bf I}_N &{\bf I}_N &{\bf I}_N &{\bf I}_N\\[0.3em]
\end{bmatrix},\]
where ${\bf I}_N$ is the $N\times N$ identity matrix and ${\boldsymbol\bbone_N}=\displaystyle\underbrace{(1\  1 \dots \ 1)}_{N\text{ times}}$, our problem is written as a standard LP problem: 

\begin{equation}\label{LPglobal}
\begin{tabular}{ r l }
$\min$ & ${\bf c}^\textup{T} {\bf x}$, \\ 
\text{subject to} & ${\bf Mx}={\bf b}$, \\ 
 & ${\bf x}\geq 0.$
\end{tabular}
\end{equation}
The result of the algorithm is a vector ${\bf x}^*:=\arg\min {\bf c}^\textup{T} {\bf x}$ whose elements $x^*_{jk}$ represent how much mass moves from node $j$ to node $k$ employing the minimum-cost mass rearrangement. 
\begin{remark}\label{dimensionGLP}
The dimension of the LP problem \eqref{LPglobal} is given by the dimensions of the matrix/vectors involved:
$$
\dimvett {\bf x}= \dimvett {\bf c} = N^2, \qquad \dimvett {\bf b}=2N, \qquad \dimvett {\bf M}=2N^3.
$$
\end{remark}
\begin{remark}
Hereafter we will refer to problem \eqref{LPglobal} as \emph{global}, in order to stress the fact that it is possible to move the mass from and to any node of the graph.
\end{remark}

\section{Coupling DMD and optimal mass transfer problem}\label{sec:coupling}

In this section we describe how we can drastically reduce the size of the LP problem \eqref{LPglobal} by using the DMD method. The resulting algorithm will allow us to study the mass transfer problem on large domains.

\subsection{Main idea}\label{sec:completeAlgorithm}
The large size of the LP problem \eqref{LPglobal}, stressed in Remark \ref{dimensionGLP}, mainly comes from the fact that the mass is allowed to be transferred from any graph node to any other graph node. However physical constraints prevent this to happening: assuming that the maximal velocity of the mass is $V_{\textup{max}}$ and denoting by $\Delta t$ the time frame from one snapshot to the following one, the maximal distance travelled is $V_{\textup{max}}\Delta t$. 

In order to reduce the size of the LP problem \eqref{LPglobal}, we restrict the set of all possible movements trajectories. The ideal time frame $\delta t$ of the snapshots would be the one which guarantees that the CFL-like condition (see, e.g., \cite[Chapter 14]{QV94})
\begin{equation}\label{CFL}
V_{\textup{max}}\delta t \leq \Delta x 
\end{equation} 
holds true. Indeed, under condition \eqref{CFL} the \emph{mass is allowed to move only towards the adjacent cells or not move at all}.

Consider now a generic set of mass distributions $\mathcal{X}=\{{\bf m}(t_0),\ldots,{\bf m}(t_{\nt})\}$, where $t_{j}=j\Delta t$, $j=0,\dots,\nt$, and $\Delta t$ is the time frame of the snapshots of the set $\mathcal{X}$. \emph{A priori}, the time step $\Delta t$ does not necessarily satisfy condition \eqref{CFL}. In particular, a large distance in time between the snapshots means that there are not enough information for reducing the set of possible movements.
However, with the DMD we are able to reconstruct the state of the system at any time instance, even if it is not provided in the original dataset.  
The coupling of DMD and of optimal mass transfer problem is done in this order:
\begin{enumerate}[label=\emph{\roman*)}]
\item We first fix $\delta t<\Delta t$ such that it satisfies condition \eqref{CFL} and then we reconstruct the solution for each $\widetilde{t}_{j}=j\delta t$ via DMD, using Algorithm \ref{Alg_DMD}. The new set of snapshots we work with is 
$\widetilde{\mathcal{X}}=\{{\bf m}(\widetilde{t}_{0}),\ldots,{\bf m}(\widetilde{t}_{\widetilde{\nt}})\}$, where ${\bf m}(\widetilde{t}_{j})$ is computed from \eqref{eq:omegaj}, $j=0,\dots,\widetilde{\nt}$. We observe that $\widetilde{\nt}>\nt$.
\item We recover the flows from the new set $\widetilde{\mathcal{X}}$ by means of an approximation of Wasserstein distance similar to the one done in Section \ref{sec:hitchcock}, but with reduced size, as described in detail below. Note that, since $\widetilde{\nt}> \nt$, we have to solve more LP problems with respect to the global approach, but, despite this, we will get advantages by using this new approach.
\end{enumerate}

Let us denote by $d$ the maximum number of neighbors per cell. 
The new unknown $\widetilde{x}_{jk}$, corresponding to the mass to be moved from node $j$ to node $k$, is defined only if $j$ and $k$ are adjacent or $j$ is equal to $k$. Analogously we define the cost function $\widetilde{c}_{jk}$.
We denote by $\bf \widetilde{x}$ the vector of the unknowns and by $\bf \widetilde{c}$ the vector associated to the cost function. We introduce the vector ${\bf s}$ whose components are indexes of nodes. The first components are the indexes of the nodes adjacent to the node 1, then those adjacent to the node 2 and so until the node $N$. Vectors ${\bf\widetilde{x}}$ and ${\bf\widetilde{c}}$ are ordered similarly to ${\bf s}$.
Since the mass can move only towards a maximum of $d$ nodes, the dimensions of $\bf\widetilde{x}$ and $\bf\widetilde{c}$ are lower or equal than $d N$.
We introduce the matrix 
\[
\bf \widetilde{M} = 
\begin{bmatrix} 
\bbone_{r_{1}} & 0 & 0 & \dots & 0\\[0.3em]
0 & \bbone_{r_{2}}& 0 & \dots & 0\\[0.3em]
0 & 0 & \bbone_{r_{3}} &  \dots & 0\\[0.3em]
\vdots & \vdots & \vdots & \ddots & \vdots\\[0.3em]
0 & 0 & 0 & \dots & \bbone_{r_{N}}\\[0.3em]
&&\bf\widetilde{D}
\end{bmatrix},
\]
where ${\bf \bbone}_{r_{i}}=\displaystyle\underbrace{(1\  1 \dots \ 1)}_{r_{i}\text{ times}}$, with $r_{i}$ the number of adjacent nodes for each node $i$, $r_{i}\leq d$, 
and ${\bf\widetilde{D}}_{ij}=1$ if the $j$-th element of ${\bf s}$ is equal to $i$, otherwise ${\bf\widetilde{D}}_{ij}=0$.
Defining the vector ${\bf \widetilde{b}} = (m^\init_1,\dots,m^\init_N,m^\fin_1,\dots,m^\fin_N)^\textup{T}$ the LP problem becomes
\begin{equation}\label{LPlocal}
\begin{tabular}{ r l }
$\min$ & ${\bf\widetilde c}^\textup{T} \bf\widetilde x$, \\ 
\text{subject to} & ${\bf\widetilde M\widetilde x}=\bf \widetilde b$ \\ 
 & ${\bf \widetilde x}\geq 0.$
\end{tabular}
\end{equation}
\begin{remark}\label{dimension_global_pb}
The dimension of the LP problem \eqref{LPlocal} is given by the dimensions of the matrix/vectors involved:
$$
\dimvett {\bf\widetilde x}= \dimvett {\bf\widetilde c} \leq dN, \qquad \dimvett {\bf\widetilde b}=2N, \qquad \dimvett {\bf\widetilde M} \leq 2dN^2.
$$
\end{remark}
\begin{remark}\label{def:global_pb}
Hereafter we will refer to problem \eqref{LPlocal} as \emph{local}, in order to stress the fact that it is possible to move the mass from and to adjacent nodes of the graph only.
\end{remark}

\subsection{A toy example for the complete algorithm: the advection equation}\label{sec:toyexamplefull}
In this test we propose an example for the complete algorithm described in Section \ref{sec:completeAlgorithm}.
Let us consider the advection equation: 
\begin{equation}\label{advection}
\begin{cases}
\partial_{t}u({\bf x},t) + {\bf v}\cdot\nabla u({\bf x},t) =0 &\quad {\bf x}\in \Omega, t\in[0,T]\\
u({\bf x},t) = 0 &\quad {\bf x}\in \partial\Omega\\
u({\bf x},0)=u_0({\bf x})  &\quad\, {\bf x}\in\Omega,
\end{cases}
\end{equation}
with ${\bf x}=(x_1,x_2)$, $ \Omega = [-2,2]\times[-2,2]$, $u_0({\bf x})=\max(0.5-x_1^2-x_2^2,0)$ and constant velocity ${\bf v}=(0.5,0.5)$. 
It is well known that the analytical solution of \eqref{advection} is $u({\bf x},t)=u_0({\bf x}-{\bf v}t)$ provided that we set $T$ small enough to have inactive boundary conditions, as it is the case here. Hereafter we denote by \emph{reference solution} the analytical solution $u$ of \eqref{advection}.
\begin{figure}[h!]
\centering
\includegraphics[scale=0.14]{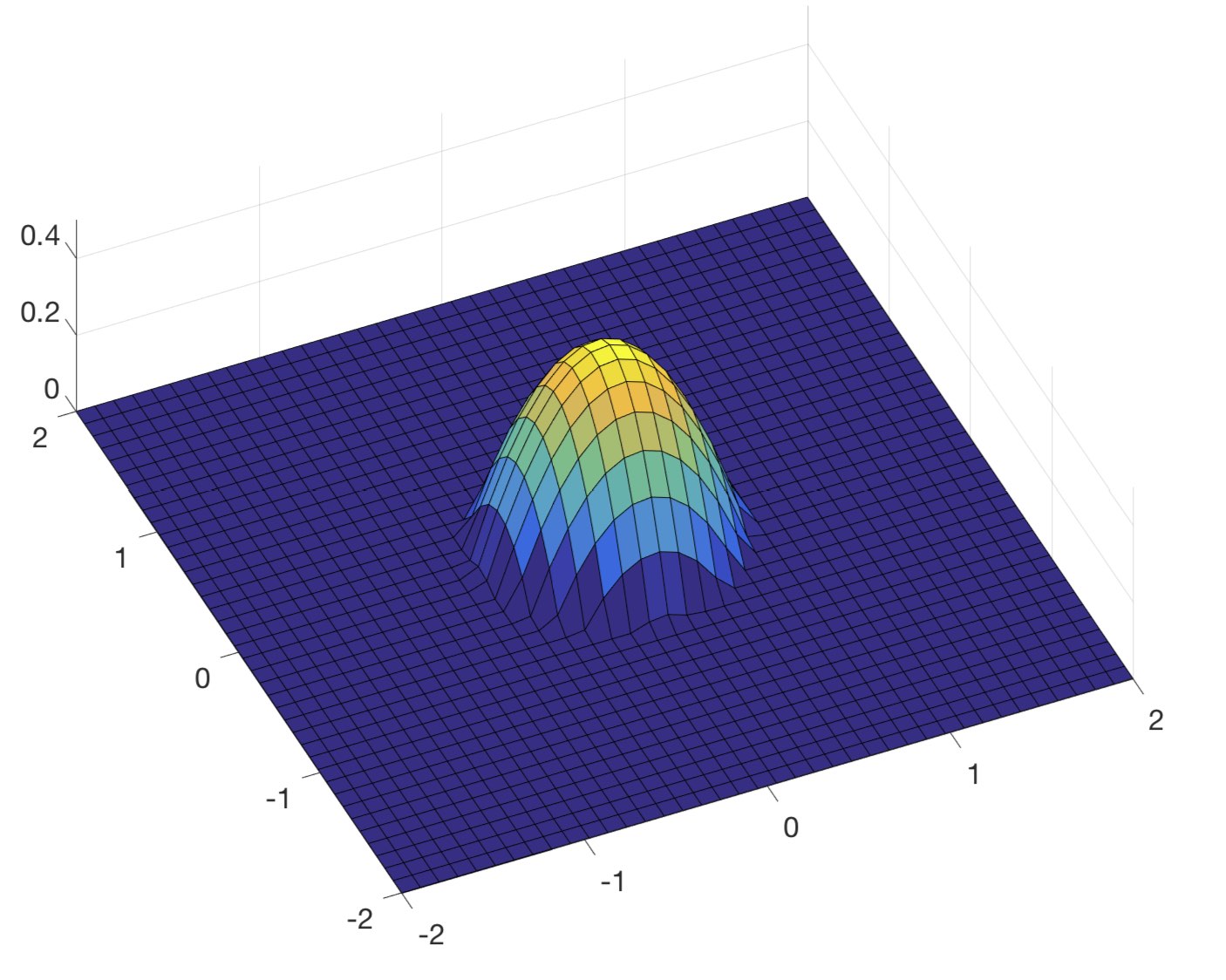}
\includegraphics[scale=0.14]{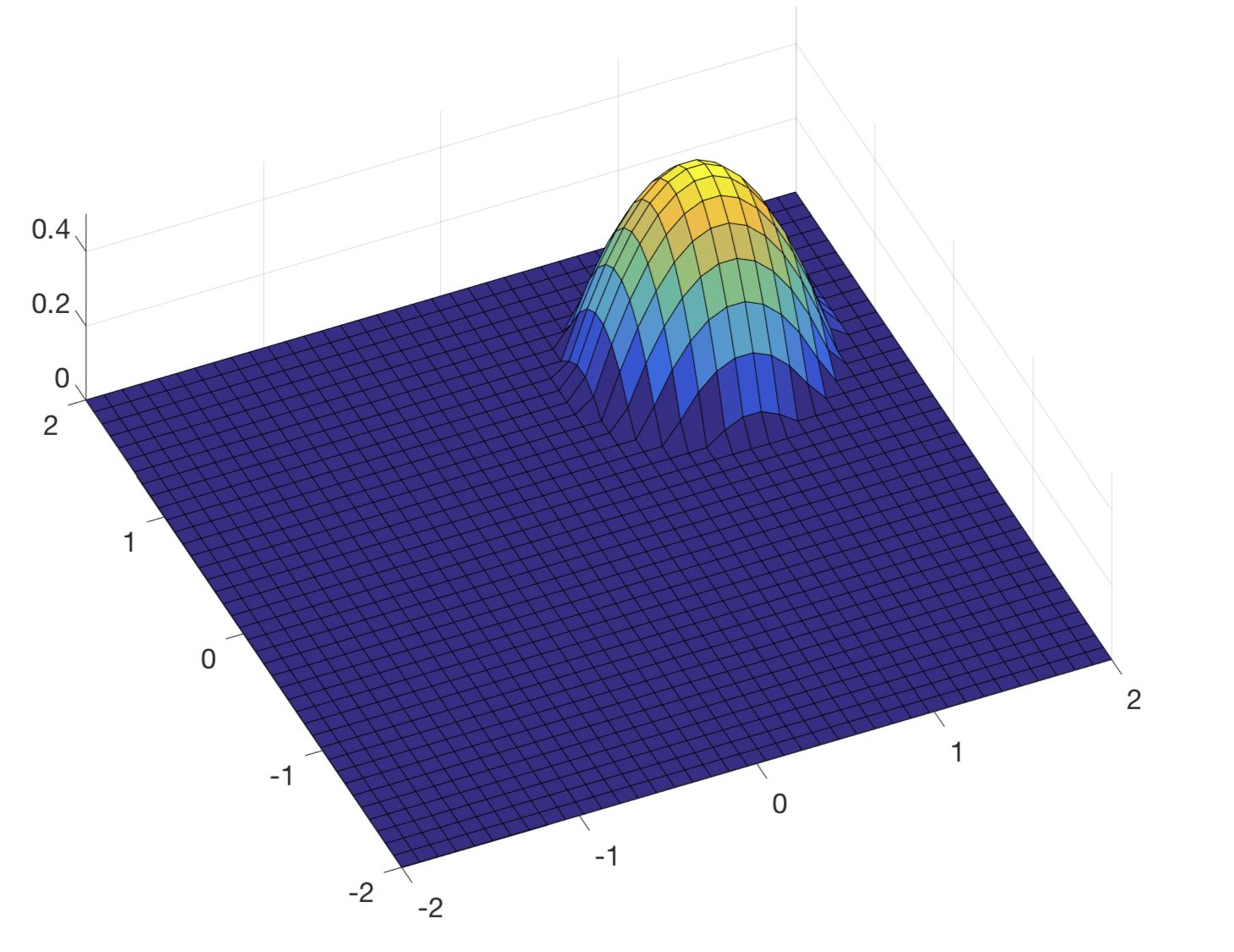}
\caption{Reference solution of equation \eqref{advection} at time $t=0$ on the left and time $t=2$ on the right.} 
\label{fig:trasportoPlot}
\end{figure}

Starting from some snapshots of the solution, we aim at reconstructing the velocity field driving the dynamics, i.e.\ the vector ${\bf v}=(0.5,0.5)$. In doing this, we will also compare the global and the local problem in terms of CPU time, see \eqref{LPglobal} and \eqref{LPlocal} respectively.

We choose a time step $\delta t$, and we collect snapshots with a larger temporal step size $\Delta t = \kappa\delta t$ with $\kappa>1$ and we reconstruct the solution via the DMD method. 
We note that here the snapshots are computed from the analytical solution of \eqref{advection}.
In particular, we work on a grid $40\times40$ and we choose $T=2$, $\Delta x = 0.1$, $\delta t = 0.05$ and $\Delta t = 2\,\delta t$. We observe that, since $V_{\textup{max}}=\norm{\bf v}=\sqrt{2}/2$, this choice of $\delta t$ fullfills the condition \eqref{CFL}. The number of snapshots is 40 and the rank $r$ in Algorithm \ref{Alg_DMD} for the DMD reconstruction is 20.
Moreover we identify the nodes of the graph $\mathcal{G}$ of the numerical approximation of Wasserstein distance with the cells of the grid.

\begin{figure}[h!]
\centering
\includegraphics[scale=1.1]{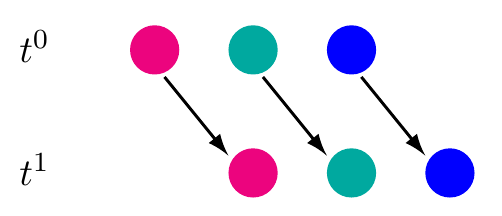}
\includegraphics[scale=1.1]{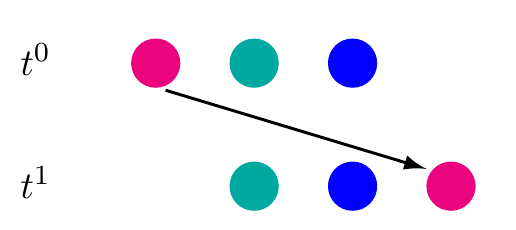}
\caption{Three small movements versus one large movement.}
\label{fig:uniqueness}
\end{figure}

\paragraph{Choice of the cost function} Since the cost function is used to figure out the ``price to pay'' for moving the mass from a node of the graph to another one, the most intuitive definition of $c$ is the Euclidean distance in $\R^2$. This choice was proved to be unsuitable for the global problem, see \cite{balzotti2018IFAC}. Indeed, since the global algorithm allows any movements between the nodes of the graph, using the Euclidean distance for the function $c$ we loose the uniqueness of the optimal transfer map ${T}^*$. To see this let us assume that we have to move three unit masses one to the right. In the picture on the left of Figure \ref{fig:uniqueness} we move the three unit masses of one to the right while in the picture on the right we move only the first mass of three to the right. The Wasserstein distance between the two configuration is clearly the same and equal to three.

The solution proposed in \cite{balzotti2018IFAC} to fix this issue was to force the minimization algorithm to select primarly the small movements penalizing the large ones. To get this, the cost function was defined as 
\begin{equation}
c(\boldsymbol{\xi}_1,\boldsymbol{\xi}_2)=\|\boldsymbol{\xi}_1-\boldsymbol{\xi}_2\|_{\mathbb R^2}^{1+\varepsilon}, 
\label{eq:costoOld}
\end{equation}
where $\boldsymbol{\xi}_1$ and $\boldsymbol{\xi}_2$ are the coordinates of the nodes of the graph and $\varepsilon>0$ is a small parameter which accounts for the penalization.\\

\begin{figure}[h!]
\centering
\subfloat[][{Local approach with reference solution.}]
{\includegraphics[scale=0.186]{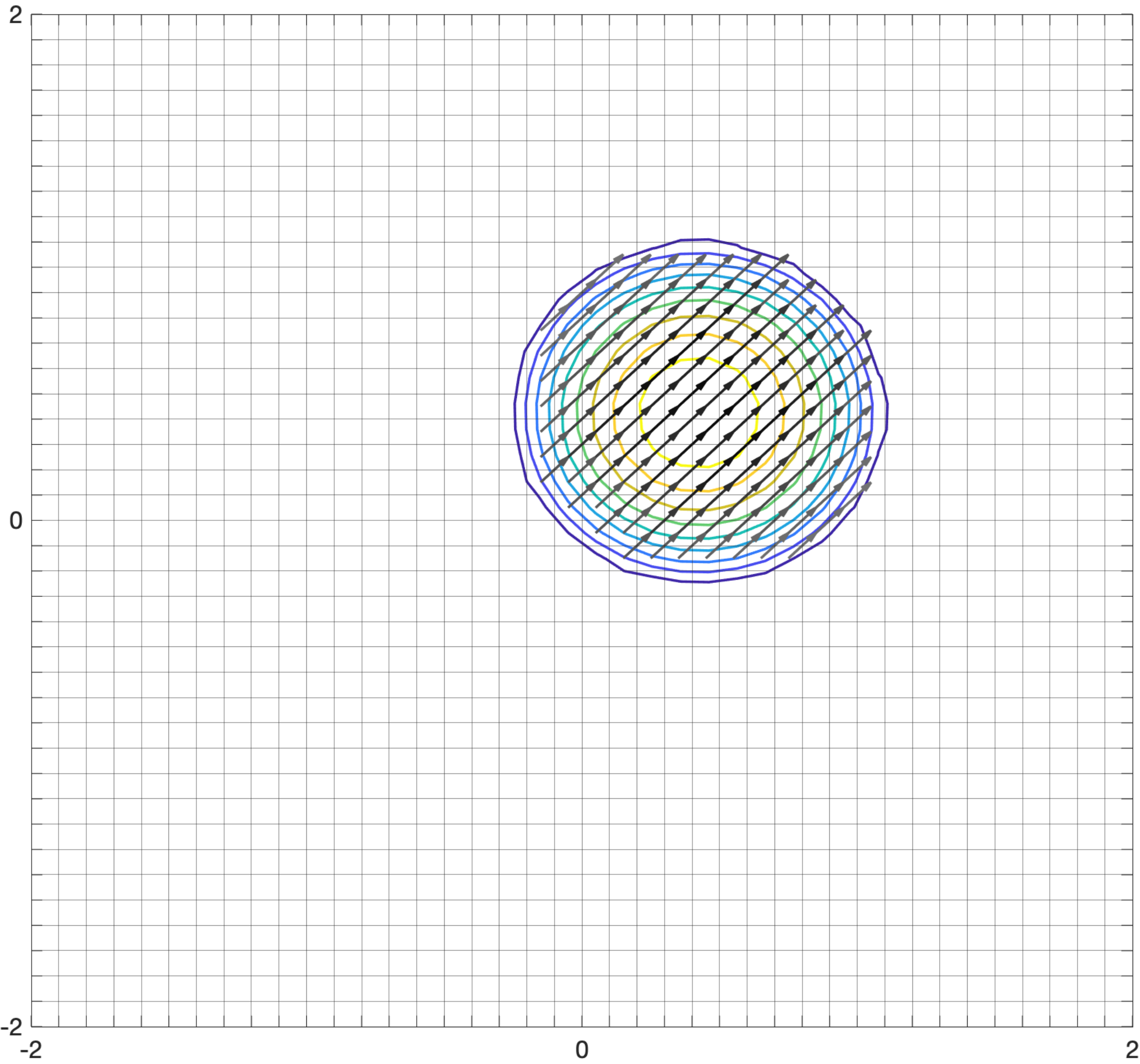}}\quad
\subfloat[][{Global approach with reference solution and Euclidean distance.}]
{\includegraphics[scale=0.186]{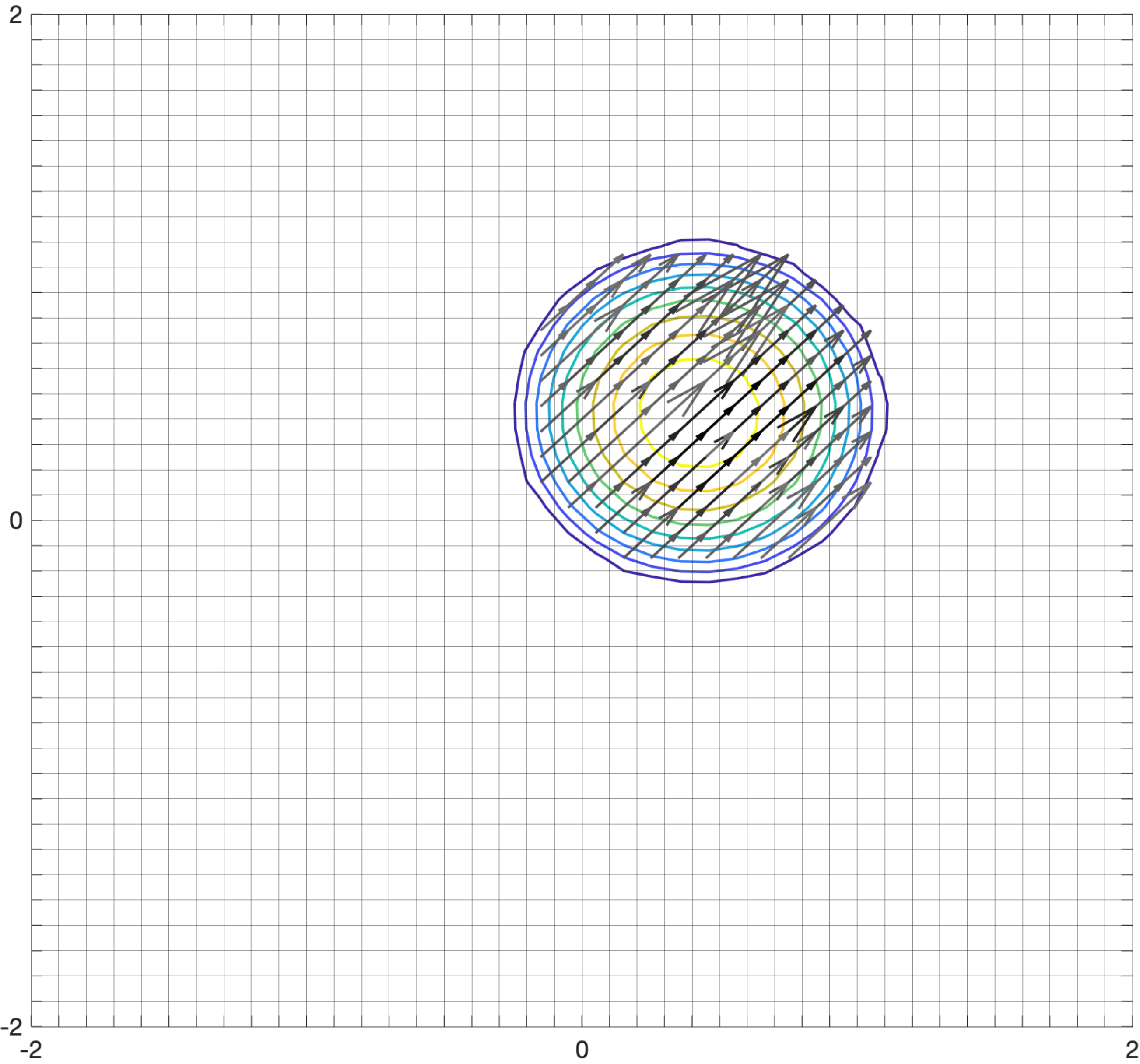}}\quad
\subfloat[][{Global approach with reference solution and cost defined in \eqref{eq:costoOld}.}]
{\includegraphics[scale=0.186]{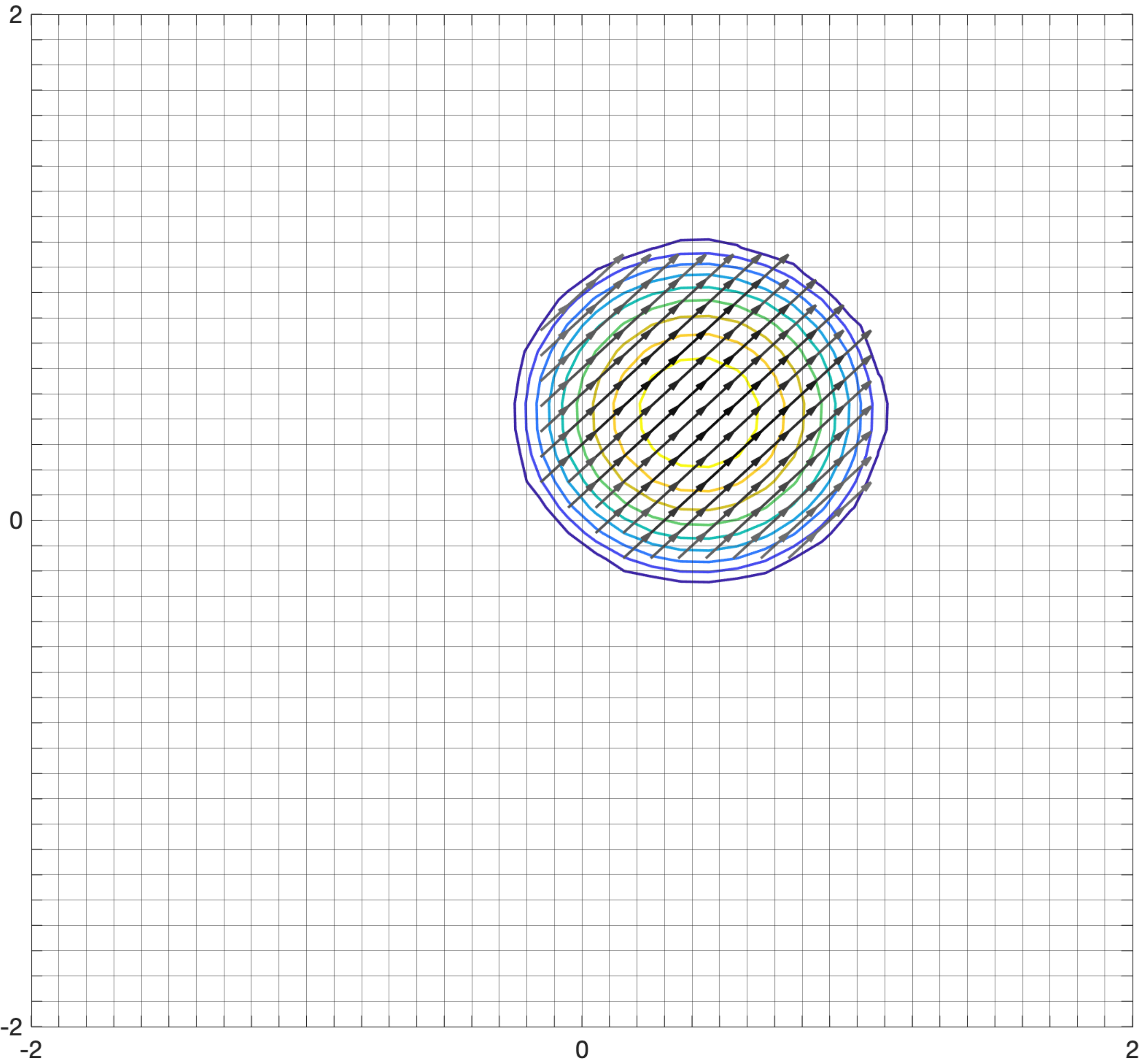}}\\
\subfloat[][{Local approach with DMD solution.}]
{\includegraphics[scale=0.186]{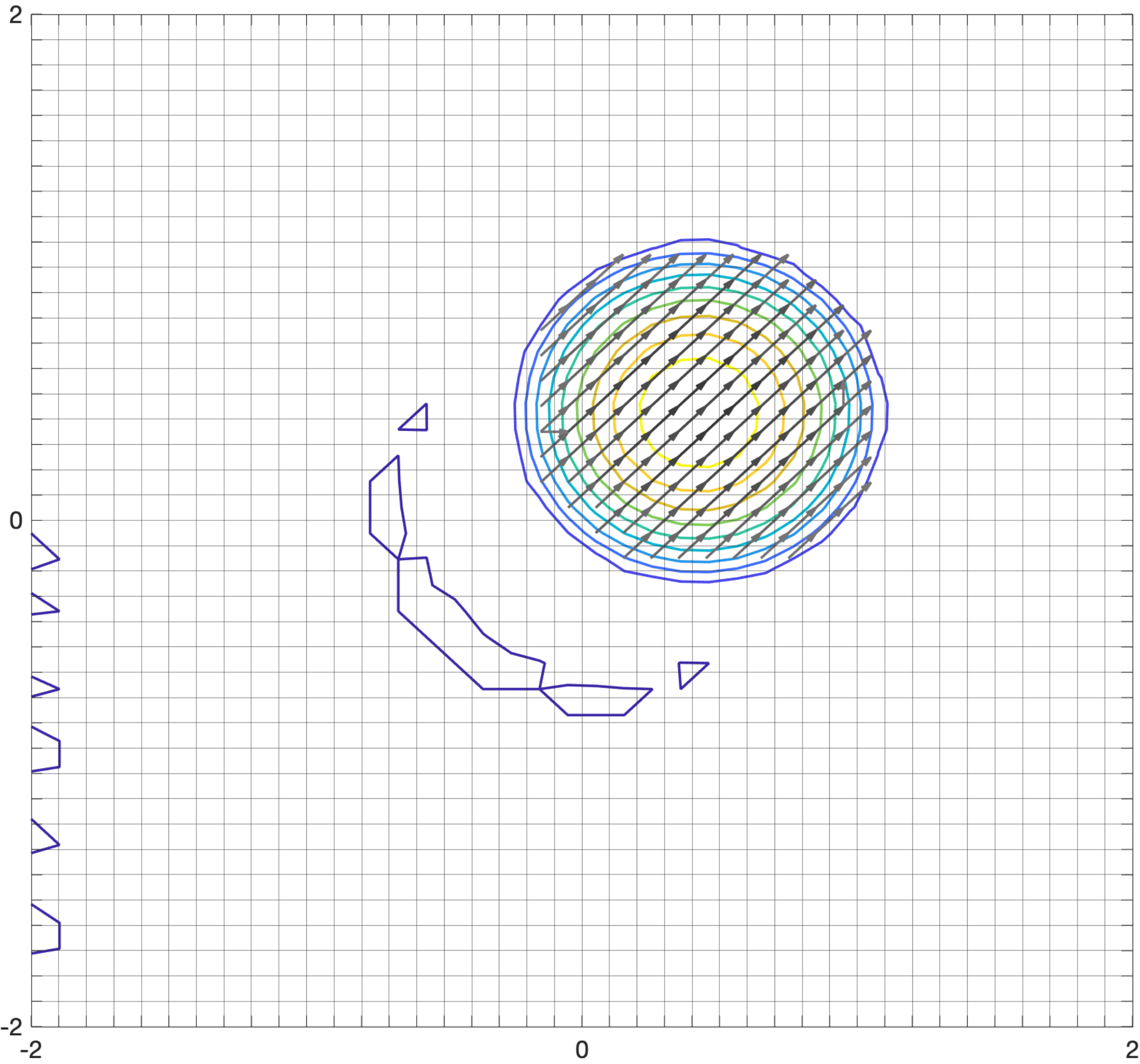}}\quad
\subfloat[][{Global approach with DMD solution and Euclidean distance.}]
{\includegraphics[scale=0.186]{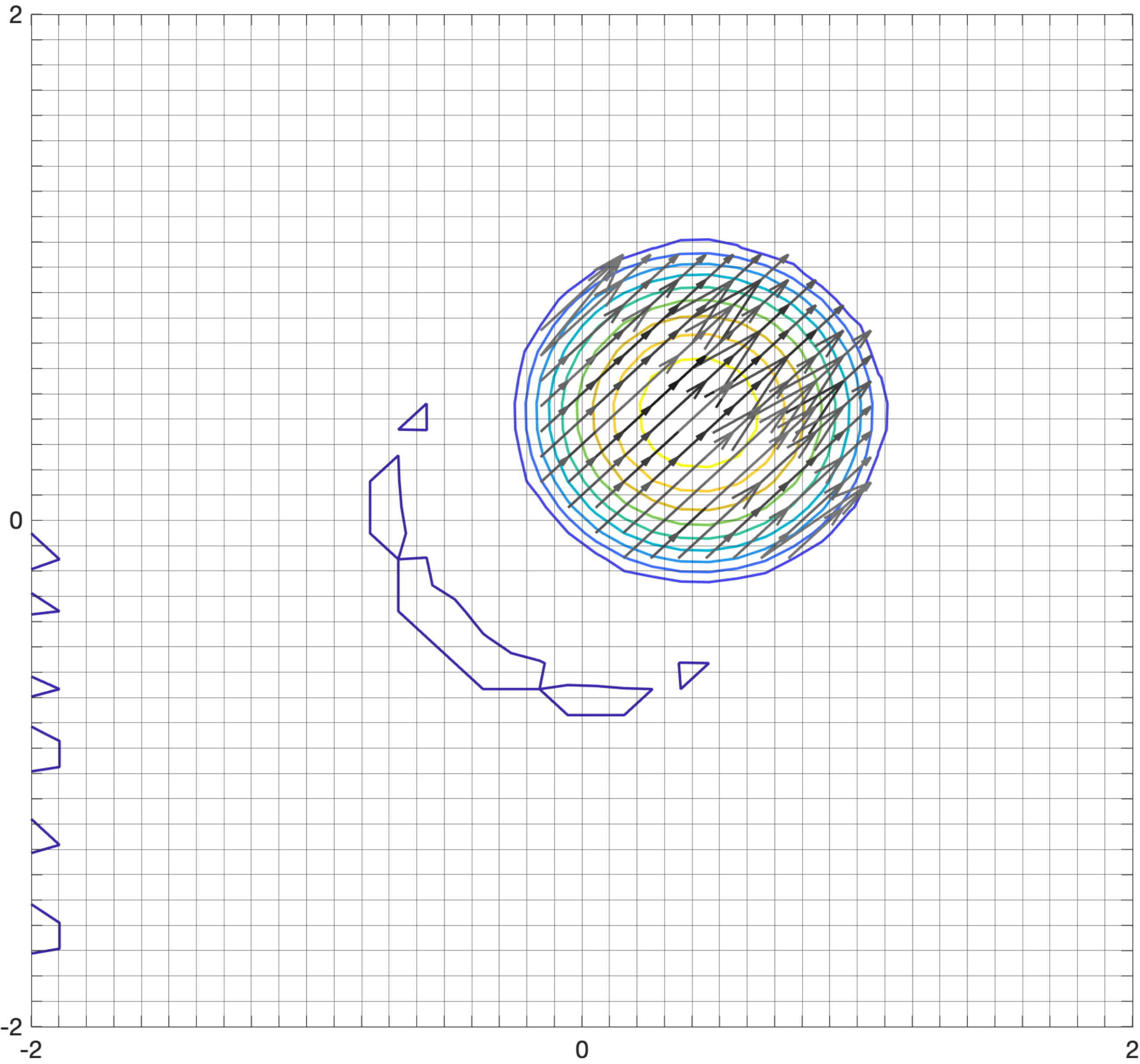}}\quad
\subfloat[][{Global approach with DMD solution and cost defined in \eqref{eq:costoOld}.}]
{\includegraphics[scale=0.186]{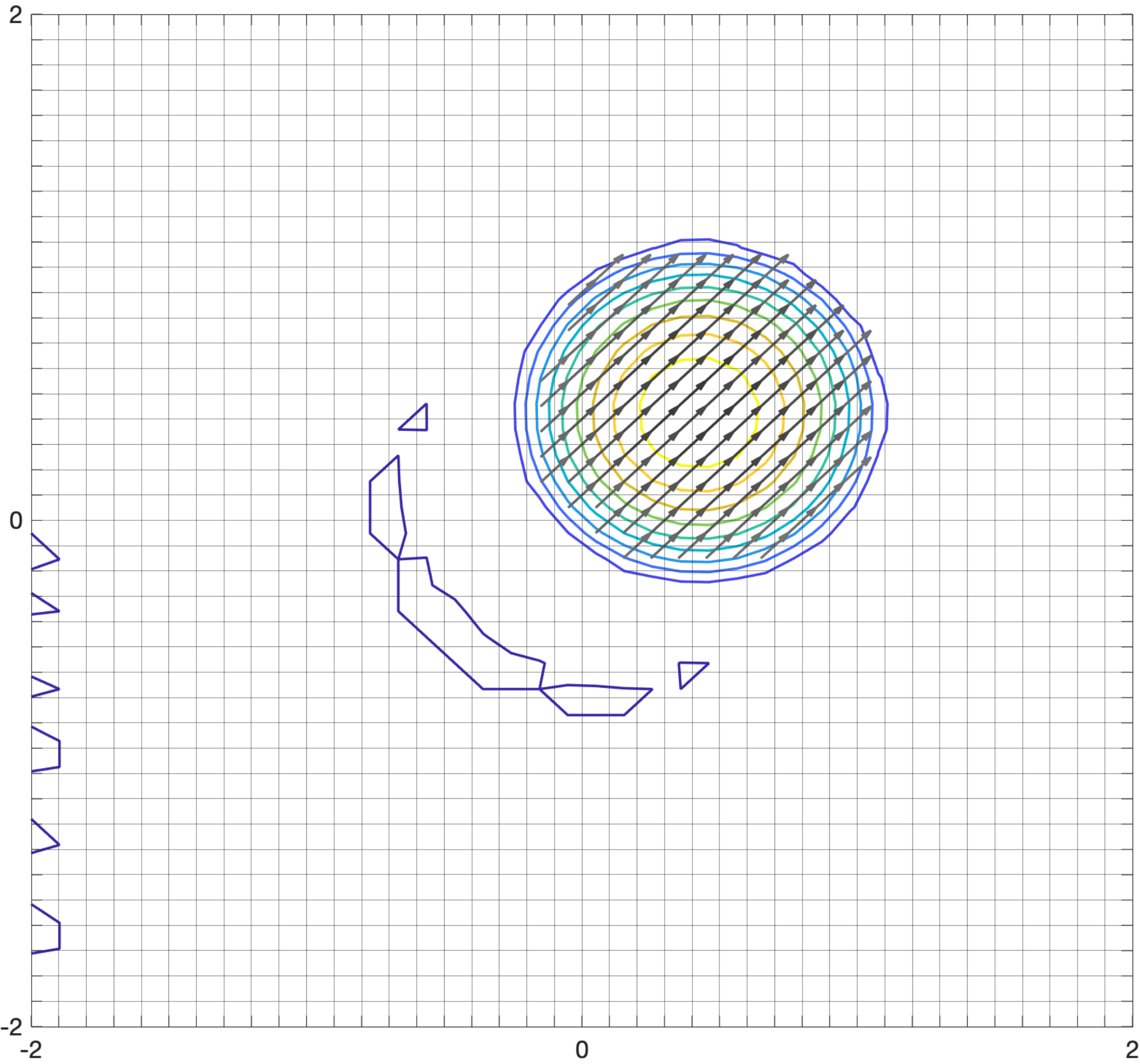}}\\
\caption{Reconstructed flows between $t_1=0.725$ and $t_2 = t_1+\delta t$ superimposed to the level sets of the solution to \eqref{advection} at $t_2$. 
Left column: local algorithm.
Central column: global algorithm with $c$ equal to the Euclidean distance.
Right column:  global algorithm with $c$ defined in \eqref{eq:costoOld}.
Top row: the reference solution to \eqref{advection} is used for computation.
Bottom row:  the DMD reconstruction of the solution to \eqref{advection} is used for computation.
}
\label{fig:trasporto}
\end{figure}

In Figure \ref{fig:trasporto} we show the level sets of the solution to \eqref{advection} together with some arrows indicating the reconstructed velocity field $\bf v$. 
More precisely, on the left column of figure we show the results obtained with the local algorithm \eqref{LPlocal}, on the center with the global algorithm with Euclidean distance and on the right with the global algorithm with the cost function defined in \eqref{eq:costoOld} (with $\varepsilon = 0.1$).
Similarly, the panels on the top show the results obtained with the reference solution whereas the panels on the bottom those obtained with the DMD solution.

The local algorithm and the global one with the correction in \eqref{eq:costoOld} are able to reconstruct the velocity field ${\bf v}=(0.5,0.5)$ with accurate approximation. The global algorithm with Euclidean distance, instead, is less precise, since the optimal flow does not correspond to the velocity field.
Moreover, the algorithm based on the DMD reconstruction of the solution introduces small oscillations in the solution to \eqref{advection}. This is expected in such hyperbolic problems since the decay of the singular values of the dataset is very slow and the initial condition is non-smooth. However, such oscillations do not have an effect in the reconstruction of the flow.

\medskip

To further validate our approach we compute the numerical error of the proposed algorithm. Since the results obtained with the global algorithm \eqref{LPglobal} with the Euclidean distance as cost function are the least accurate, we focus only on the other two approaches.
For each time step $t_{n}=n\delta t$ we define ${\bf x}^*_{E}(t_{n})$ and ${\bf\widetilde{x}}^*_{E}(t_{n})$ as the solutions to the LP problems \eqref{LPglobal} and \eqref{LPlocal} respectively at time $t_{n}$, when the known terms ${\bf b}(t_{n})$ and ${\bf\widetilde{b}}(t_{n})$ are chosen as the reference solution to \eqref{advection} at time $t_{n}$.
Analogously, ${\bf x}^*_{D}$ and ${\bf\widetilde{x}}^*_{D}$ are the solution when the known terms are obtained with the DMD method. 
The vectors ${\bf x}^*_{E}(t_{n})$, ${\bf\widetilde{x}}^*_{E}(t^{n})$, ${\bf x}^*_{D}(t^{n})$ and ${\bf\widetilde{x}}^*_{D}(t_{n})$, for $n=1,\dots,\floor{\frac{T}{\delta t}}+1$, are finally collected as the columns of the matrices ${\bf X}^*_{E}$, ${\bf\widetilde{X}}^*_{E}$,  ${\bf X}^*_{D}$ and ${\bf\widetilde{X}}^*_{D}$ respectively. 
We define the errors 
\begin{equation}\label{eq:erroreTrasporto1}
E^{G} := \frac{\norm{{\bf X}^*_{E}-{\bf X}^*_{D}}_F}{\norm{{\bf X}^*_{E}}_F},
\qquad
E^{L} := \frac{\norm{{\bf\widetilde{X}}^*_{E}-{\bf\widetilde{X}}^*_{D}}_F}{\norm{{\bf\widetilde{X}}^*_{E}}_F},
\end{equation}  
where $\norm{\cdot}_{F}$ is the Frobenius norm.
In Table \ref{tab:trasporto} we compare the errors defined in \eqref{eq:erroreTrasporto1} obtained from the simulations on a grid with $N\times N$ nodes, for $N=20,\,30$ and $40$. As we can see from the table, increasing the number of nodes we reduce the space step $\Delta x = 4/N$ and the time step $\dt=\Delta x/2$ and
we obtain the decrease of the error between the reference solution and the DMD reconstruction. Moreover, the error obtained with the local algorithm is significantly smaller than the one obtained with the global approach.

\begin{table}[tbhp]
{\footnotesize
\caption{Comparison of the errors defined in \eqref{eq:erroreTrasporto1}.}\label{tab:trasporto}
\begin{center}
\begin{tabular}{|c|l|l|c|l|}\hline
\multicolumn{1}{|c|}{$N$} & \multicolumn{1}{c|}{$\Delta x$} &\multicolumn{1}{c|}{$\Delta t$} & \multicolumn{1}{c|}{$E^G$} & \multicolumn{1}{c|}{$E^L$}\\ \hline
20 & 0.20 & 0.100 & 0.185 & 0.030\\ 
30 & 0.13 & 0.067 & 0.159 & 0.028\\ 
40 & 0.10 & 0.050 & 0.122 & 0.020\\ \hline
\end{tabular}
\end{center}
}
\end{table}

Finally, in Table \ref{tab:trasportoCosto} we compare the computational time between the global approach \eqref{LPglobal}, with the cost function as in \eqref{eq:costoOld}, and the local approach \eqref{LPlocal} with respect to the nodes of the graph. We observe that the local algorithm is always faster than the global one. The difference between the two approaches becomes more relevant when we refine the grid by increasing the number of nodes, and thus the number of time steps. Specifically, for a grid $40\times 40$, the local algorithm required a few seconds whereas the global one more than three hours.

\begin{table}[tbhp]
{\footnotesize
\caption{Computational time.}\label{tab:trasportoCosto}
\begin{center}
\begin{tabular}{|c|r|r|r|r|} \hline
$N$ & \bf Global Exact & \bf Global DMD & \bf Local Exact & \bf Local DMD \\ \hline
20 & 18.10 s &18.89 s & 0.28 s &0.36 s\\ 
30 & 11 min & 11 min & 0.93 s & 1.41 s\\ 
40 & 3 h 6 min & 3 h 7 min & 2.40 s & 4.85 s\\  \hline
\end{tabular}
\end{center}
}
\end{table}

\section{Application to real mobile phone data}\label{sec:TIM}
In this section we focus on a specific application of the proposed approach. The real dataset gives information about the spatial distribution of people in a large populated area. The goal is to understand the travel flows of people, focusing in particular on recurring patterns and daily flows of commuters.

\subsection{Dataset}
The Italian telecommunication company TIM provides estimates of mobile phones presence in a given area in raster form: the area under analysis is split into a number of elementary territory units (ETUs) of the same size (about 130 $\times$ 140 m$^2$ in urban areas). The estimation algorithm does not singularly recognize users and does not track them using GPS. It simply counts the number of phone attached to network nodes and, knowing the location and radio coverage of the nodes, estimates the number of TIM users within each ETU at any time. TIM has now a market share of 30\% with about 29.4 million mobile lines in Italy (AGCOM, Osservatorio sulle comunicazioni 2/2017).

\begin{figure}[h!]
\begin{center}
\includegraphics[width=8cm]{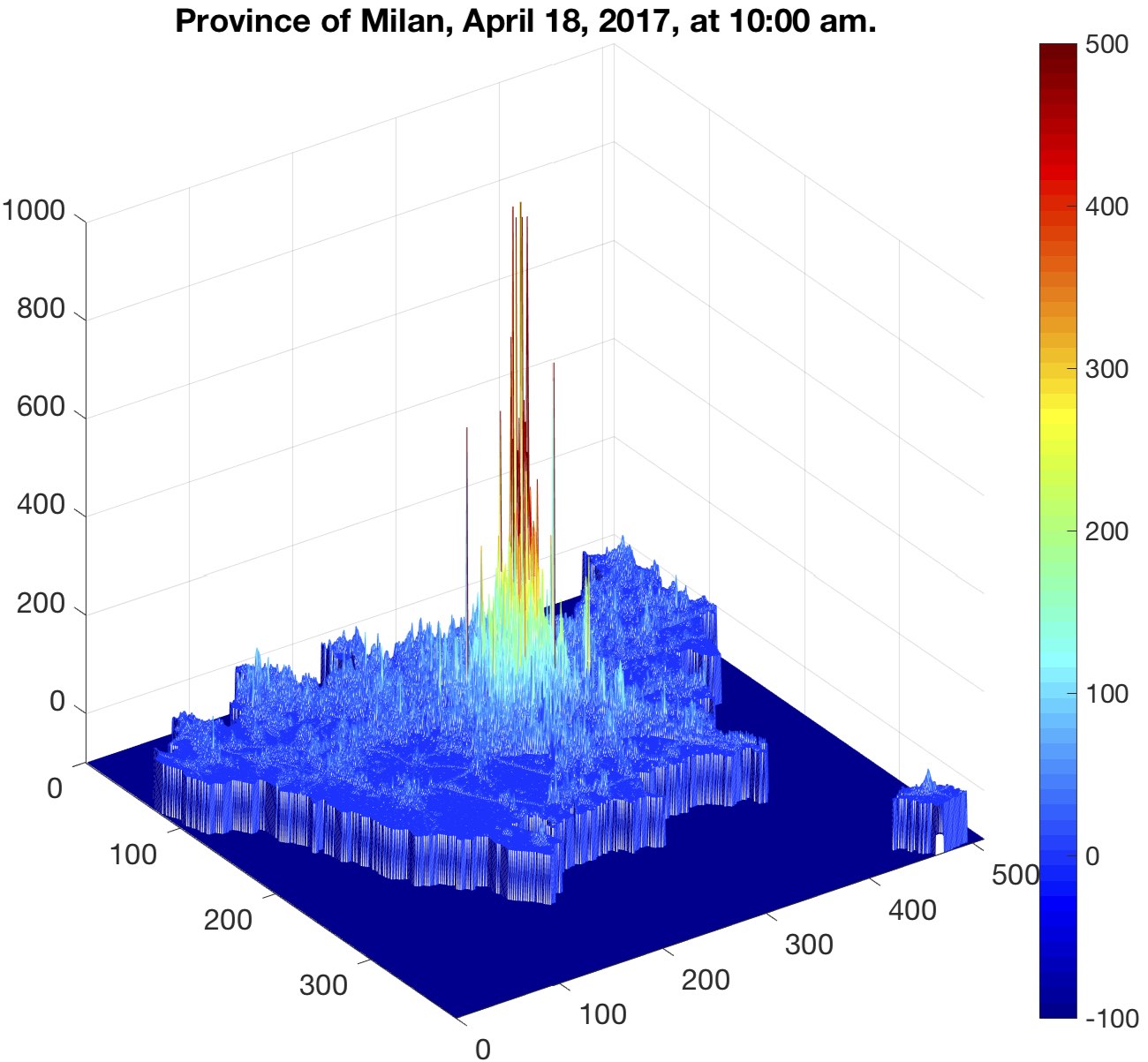}    
\caption{3D-plot of the number of TIM users in each ETU of Milan's province on April 18, 2017.} 
\label{fig:presenze3D}
\end{center}
\end{figure}

\begin{figure}[h!]
\begin{center}
\includegraphics[scale=.27 ]{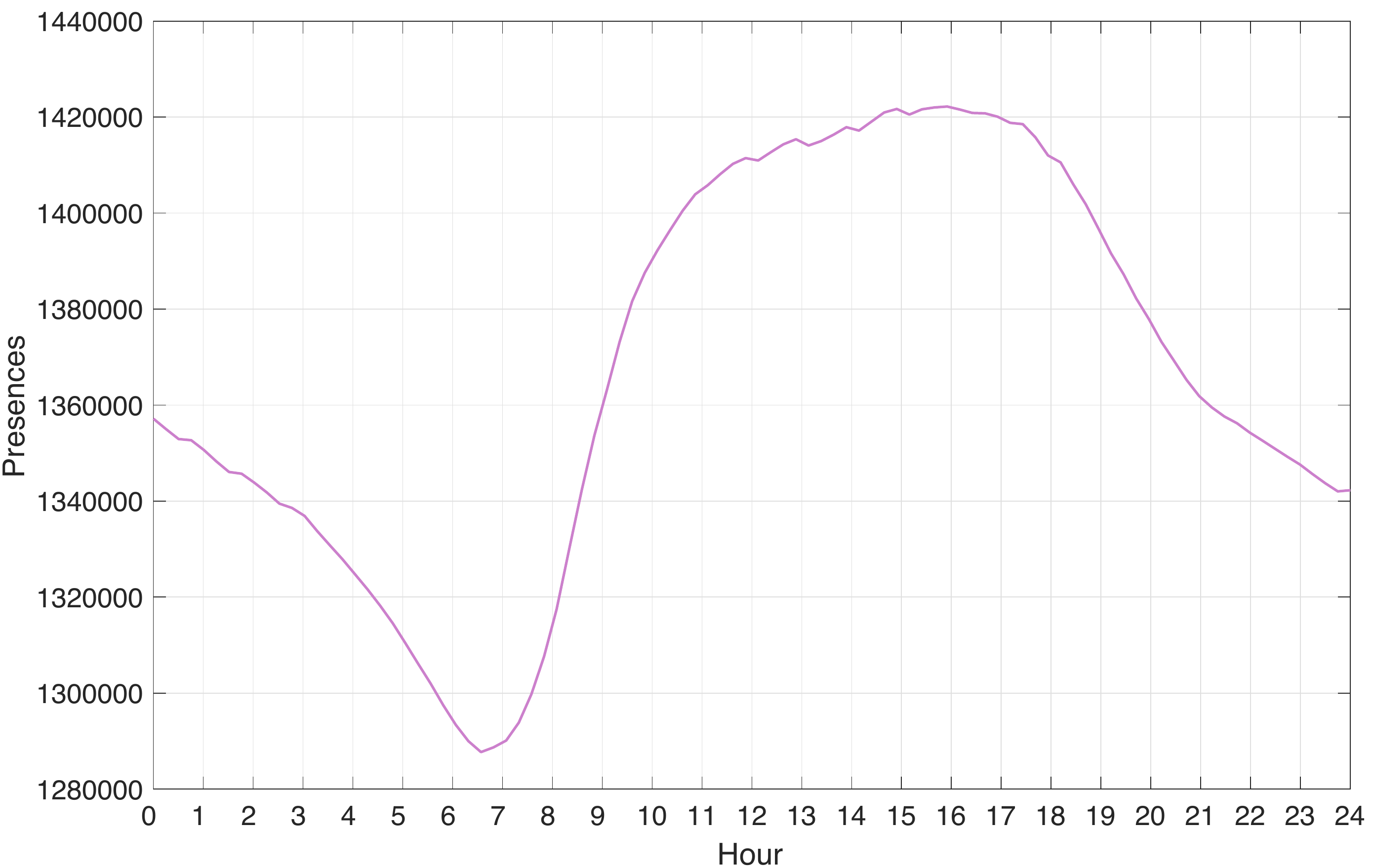}    
\includegraphics[scale = .27]{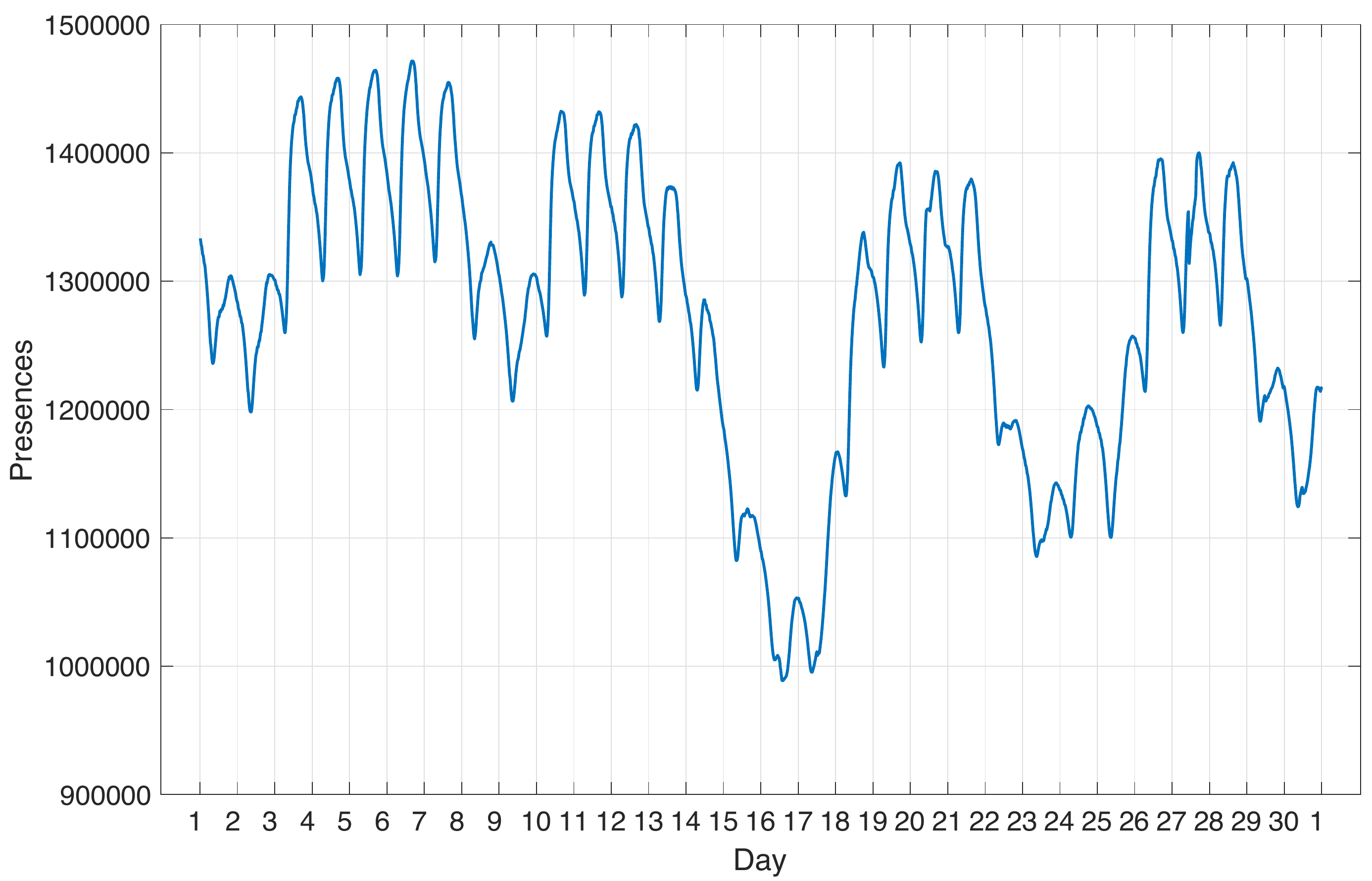}  
\caption{Trend of presences in the province of Milan during a typical working day (left). Trend of presences in the province of Milan during April 2017 (right).} 
\label{fig:presenze}
\end{center}
\end{figure}

The data we considered refer to the area of the province of Milan (Italy), which is divided in 198,779 ETUs, distributed in a rectangular grid  389 $\times$ 511. Data span six months (February, March and April 2016 and 2017). The entire period is divided into time intervals of 15 minutes, therefore we have 96 data per day per ETU in total.
In Figure \ref{fig:presenze3D} we graphically represent presence data at a fixed time. We observe that the peak of presence is located in correspondence of Milan city area.
Figure \ref{fig:presenze} shows the presences in the province of Milan in a typical working day in the left panel. The curve in the image decreases during the night, it increases in the day-time and decreases again in the evening. These variations are due to two main reasons: first, the arrival to and departure from Milan's province of visitors and commuters. Second, the fact that when a mobile phone is switched off or is not communicating for more than six hours, its localization is lost. 
The presence value that most represents the population of the province is observed around 9  pm., when an equilibrium between traveling and phone usage is reached. This value changes between working days and weekends, but it is always in the order of $1.3\times10^6$.
On the right panel of Figure \ref{fig:presenze} we show the trend of presence data during April 2017. We can observe a cyclical behavior: in the working days the number of presences in the province is significantly higher than during the weekends. It is interesting to note the presence of two low-density periods on April 15-18 and on April 22-26, 2017, determined respectively by the Easter and the long weekend for the Italy's Liberation Day holiday.
\subsection{DMD approach on TIM data}
As explained in Section \ref{sec:introDMD}, we can reconstruct the presence data in each cell at any time. We denote by ${\bf m}(t^\init)\in\R^N$ the vector containing the number of people present in the $N$ cells at a certain quarter of an hour and by ${\bf m}(t^\init+15\min)\in\R^N$ the same quantities at the consecutive time step. Applying the DMD method we are able to calculate ${\bf m}(t)\in\R^N$ for any $t\in[t^\init,t^\init+15\min]$, see \eqref{eq:omegaj}. We also note that DMD can be applied to higher dimensional dataset through randomized methods as explained in \cite{AK19}.

\begin{figure}[h!]
\centering
\includegraphics[scale = .32]{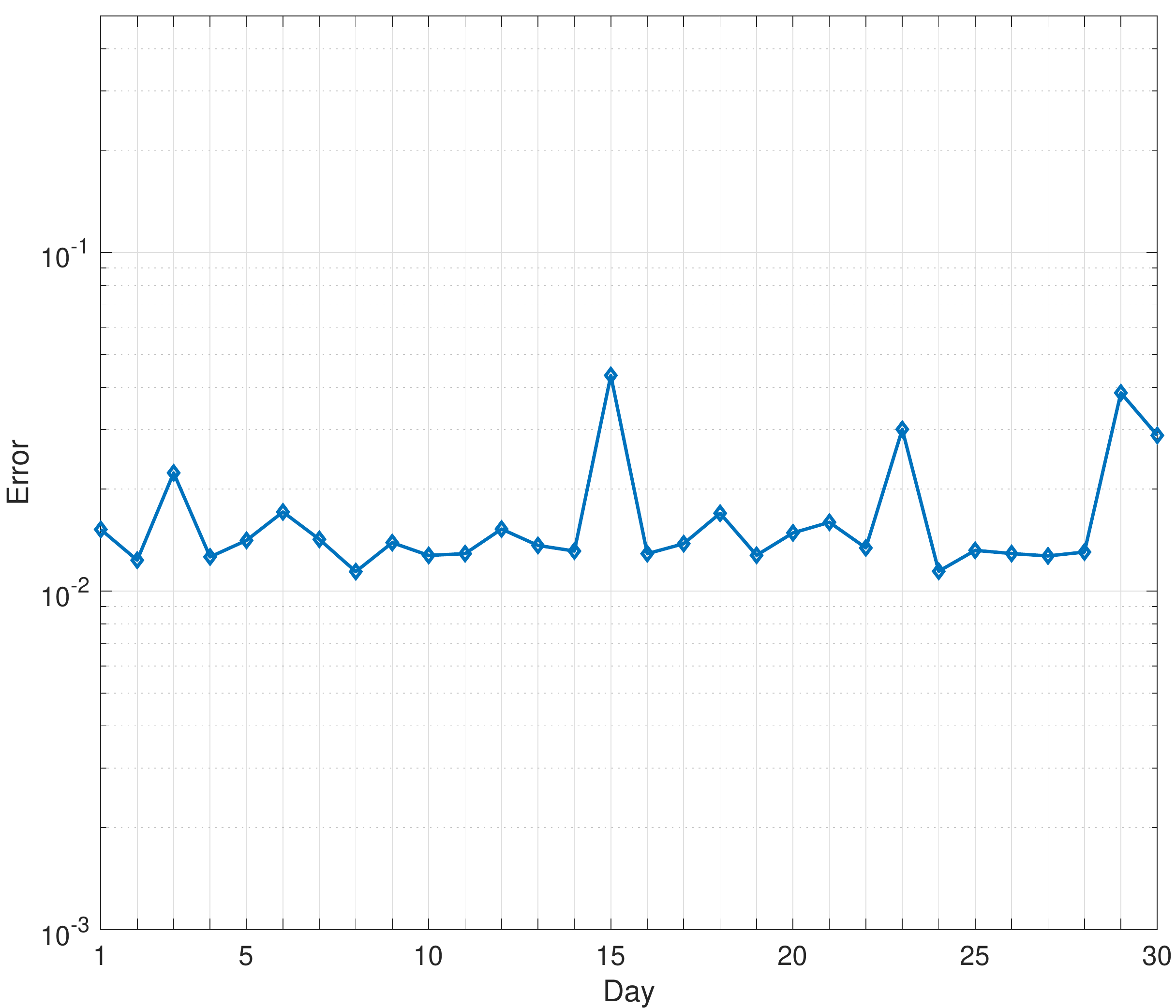}   
\caption{Daily error defined in \eqref{eq:erroreDMD} where DMD has been computed with $r=95$ in step 2 of Algorithm \ref{Alg_DMD}.}
\label{DMD_TIM}
\end{figure}

To validate the DMD approach on the TIM dataset we study the daily error  using only half of the data at our disposal in the DMD reconstruction. More precisely, we denote by $\pdata$ the matrix whose columns correspond to the real data of presences stored every $15$ minutes. 
Then, we reconstruct the data of presences every minute with the DMD algorithm, using snapshots every 30 minutes. In other words, we exploit only one column out of two of $\pdata$ to build the matrix $\pdmd$, which collects the reconstructed data every minute.
Since a day contains 96 quarters of an hour and 1440 minutes, $\pdata$ is a matrix $96\times N$, while $\pdmd$ is a matrix $1440\times N$. Finally, to compare the real data with the DMD reconstruction, we extract from $\pdmd$ the rows corresponding to the original interval of $15$ minutes into the matrix ${\bf\widetilde{P}}^{\mathrm{DMD}}$, of dimensions $96\times N$, and then we define the error as:
\begin{equation}
E = \frac{\norm{\pdata-{\bf\widetilde{P}}^{\mathrm{DMD}}}_F}{\norm{\pdata}_F}.
\label{eq:erroreDMD}
\end{equation}
In Figure \ref{DMD_TIM} we show the daily error \eqref{eq:erroreDMD} for an entire
month of data in the whole area of the province of Milan. 
The daily error is of order $10^{-2}$, which certifies the accuracy of the DMD method.

\subsection{Understanding human mobility flows}

Following the approach described in Section \ref{sec:coupling}, we define the graph $\mathcal{G}$ by exploiting the subdivision of the area of the province of Milan (Italy) into ETUs; 
we identify the $N$ nodes of the graph with the corresponding center of the ETUs, ordered from the left to the right and from the top to the bottom. 
The result is a rectangular graph $\mathcal{G}$ divided into $N_R$ rows and $N_C$ columns ($N=N_R\times N_C$). The mass $m_j(t_{n})$ is defined as the average number of presences in the node $j$ at time $t_n=n\, 15\min$.

Let us assume that $V_{\textup{max}}$ in \eqref{CFL} is equal to 50 km/h. Since the dimensions of the ETUs is around 150 m, to apply the DMD we fix the new time step $\delta t$ equal to 10 seconds. With this choice we assume that people can move only towards the eight adjacent nodes of the rectangular graph, or not move at all. We observe that the mass in the nodes on the four corners of the graph can move only towards four directions (adjacent nodes or no movement), while the mass along the boundaries can move only towards six directions. In this way, the total number of possible movements between the cells is given by 
\begin{equation}
\widetilde{N}= \underbrace{4\cdot 4}_{\text{corners}}+ \underbrace{6\cdot 2\,(N_R+N_C-4)}_{\text{boundaries}}+\underbrace{9\,(N-4-2\,(N_R+N_C-4))}_{\text{internal nodes}}<9N.
\label{eq:Ntot}
\end{equation}

The vectors $\bf\widetilde{x}$ and $\bf \widetilde{c}$ associated to the unknown  moving mass and to the cost function have length $\widetilde{N}$, while the matrix $\bf\widetilde{M}$ of the LP problem \eqref{LPlocal} has dimension $2N\times \widetilde{N}$. In Table \ref{tab:confronto}, we compare the dimensions of the vectors and the matrices of the two different LP problems: it is clear that the computational time to solve problem \eqref{LPlocal} is significantly reduced respect to problem \eqref{LPglobal}.
\begin{table}[tbhp]
{\footnotesize
\caption{Comparison between dimensions of matrices and vectors for the two methods.}\label{tab:confronto}
\begin{center}
\begin{tabular}{|c|c|c|} \hline
\bf Algoritm & \bf Vectors Dimension & \bf Matrix Dimension \\ \hline
Global & $N^2$ & $2N^3$\\
Local & $\widetilde{N}\, (<9N)$ & $2N\times\widetilde{N} \,(<18N^2)$\\ \hline
\end{tabular}
\end{center}
}
\end{table}

\paragraph{Choice of the cost function $c$}
Since the ETUs are rectangles of length $\ell_x=130$ m and $\ell_y=140$ m, we define the cost function $c_{jk}$ for the local algorithm as follows:
\begin{equation*}
c_{jk}=\begin{cases}
0 &\quad\text{if $j=k$ or $j$ and $k$ are not adjacent}\\
\ell_x &\quad\text{if $j$ and $k$ are horizontally adjacent}\\
\ell_y &\quad\text{if $j$ and $k$ are vertically adjacent}\\
\sqrt{\ell_x^2+\ell_y^2} &\quad\text{if $j$ and $k$ are diagonally adjacent.}
\end{cases}
\end{equation*}
For the global approach, we use the cost function defined in \eqref{eq:costoOld} with $\varepsilon = 0.1$.

To sum up, in order to solve the mass transfer problem for a whole day using snapshots taken every 15 minutes (real data) we have to solve 96 global LP-problems \eqref{LPglobal}, whereas with the DMD algorithm we have to solve 8640 local LP-problems \eqref{LPlocal}. As we will see in the following section, despite the larger number of LP problems, the local approach is more convenient than the global one.
\begin{remark}
The Wasserstein distance is defined between two distributions of equal mass (see \eqref{def:WassDist}). In our case the conservation of mass between two consecutive snapshots is not guaranteed. 
Let us consider a couple of snapshots with different total mass $\sum_{j} m_{j}^{0} \neq\sum_{j} m_{j}^{1}$. To correctly apply the algorithms for the identification of flows, we compute the mass in excess between the two snapshots and then we uniformly distribute it in all the nodes of the graph with lower mass. 
A more sophisticated approach could be the one suggested in \cite{piccoli2014ARMA}, where a definition of Wasserstein distance between two distributions with different mass is proposed.
\end{remark}

\subsection{Numerical results}
In this section we show the results obtained with the local algorithm \eqref{LPlocal} to study the flows of commuters and the influence of great events on human mobility.  In both cases the rank $r$ in step 2 of Algorithm \ref{Alg_DMD}, used for the construction of the DMD solution, is 95.
The flows are represented by arrows; we draw only those which correspond to the most significant movements, and we associate a darker color to the arrows corresponding to a  larger movement of people. For graphical purposes, in the following plots we aggregate 6 time steps $\delta t$ to show 1-minute  mass transport.

\subsubsection{Flows of commuters}
\begin{figure}[h!]
\centering
\includegraphics[scale=0.56]{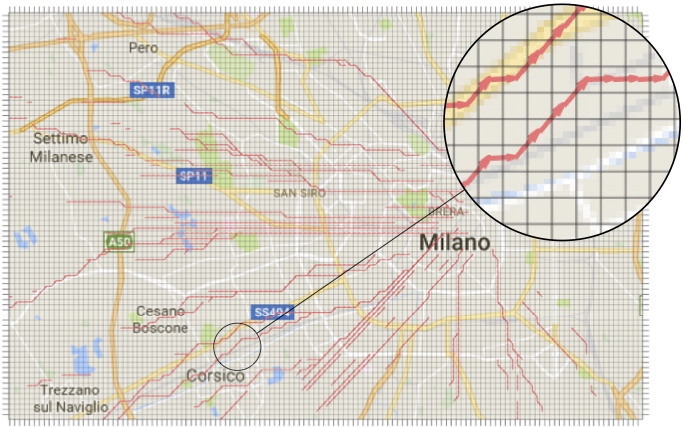}
\includegraphics[scale=0.56]{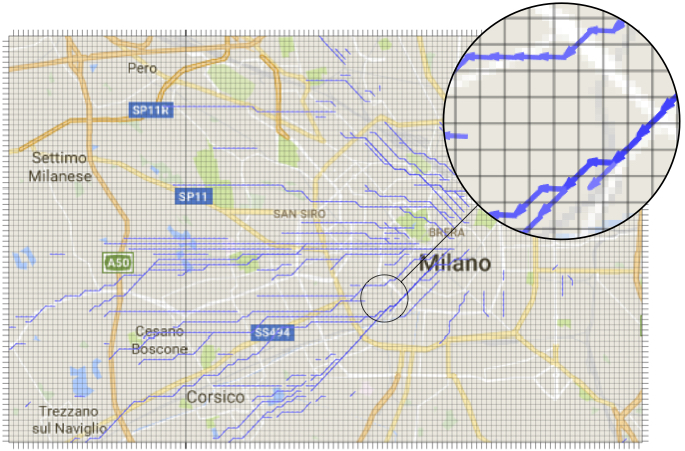}
\caption{Flows of commuters during the morning (09:00-09:01) of a generic working day (top). Flows of commuters during the evening (18:00-18:01) of a generic working day (bottom).}
\label{fig:provincia}
\end{figure}

In the following simulations we consider the area of the Province of Milan during a generic working day. Milan is one of the biggest Italian city and it attracts many workers from outside. The city of Milan is located in the right part of the Province, therefore we mainly see movements from/to the left part of the analyzed area. In the top panel of Figure \ref{fig:provincia} we show the morning flows of a working day: we clearly see that the arrows are directed towards the city of Milan. In particular, we zoom on the arrows which overlap the roads heading to Milan. In the bottom panel of Figure \ref{fig:provincia} we show the opposite phenomenon: in the evening people go away from work to come back home inside the Province of Milan.

\paragraph{CPU times}
Considering data for a 6 hours frame on an area of $144\times 240$ ETUs the local algorithm requires 144 hours of CPU time and works with 360 snapshots. The global approach \eqref{LPglobal} is not able to analyze such an area, since the matrix $\bf M$ in \eqref{LPglobal} has a computationally unmanageable dimension.

\subsubsection{Flows influenced by a great event}
\begin{figure}[h!]
\centering
\includegraphics[width=.48\columnwidth]{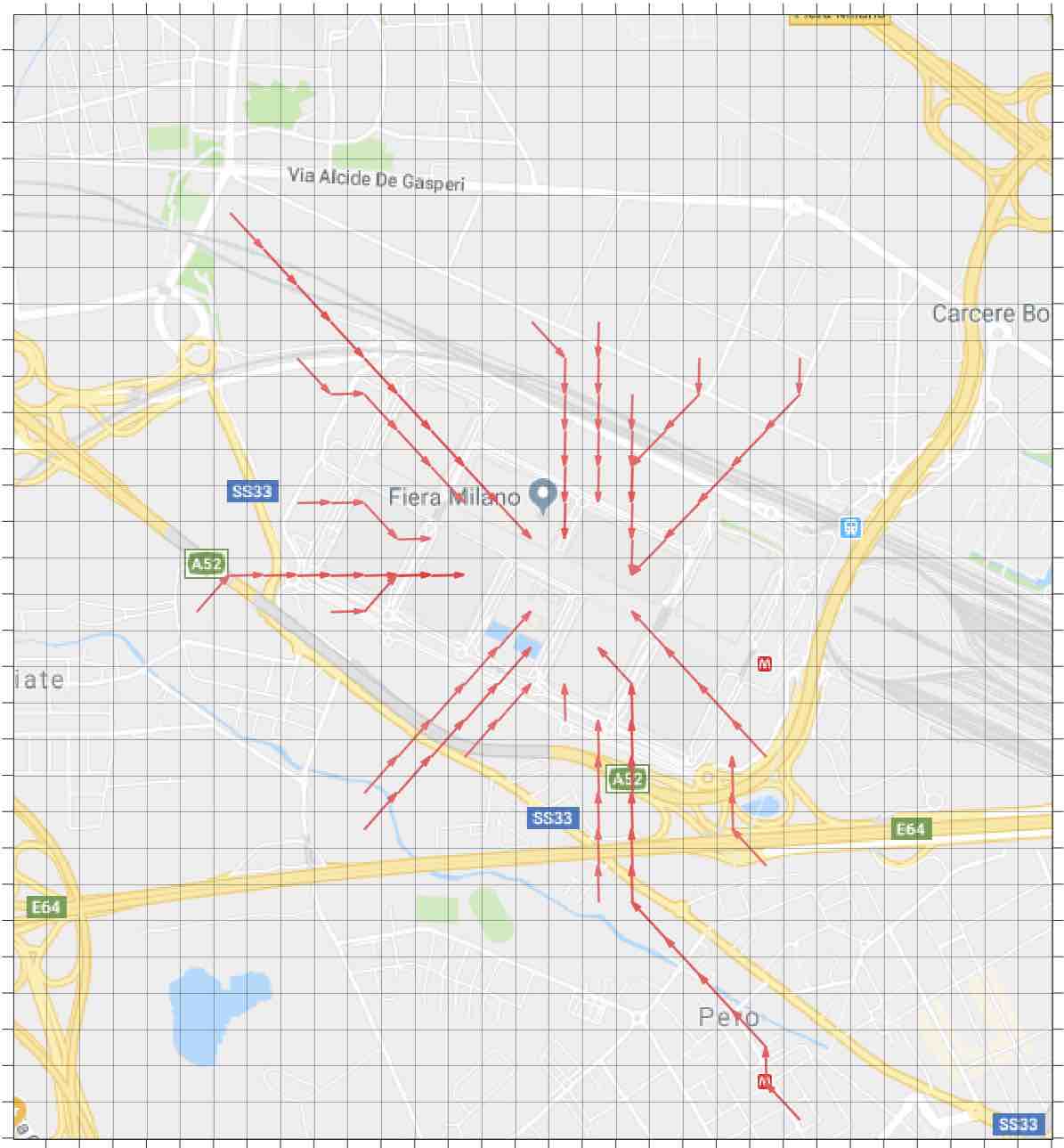}
\includegraphics[width=.48\columnwidth]{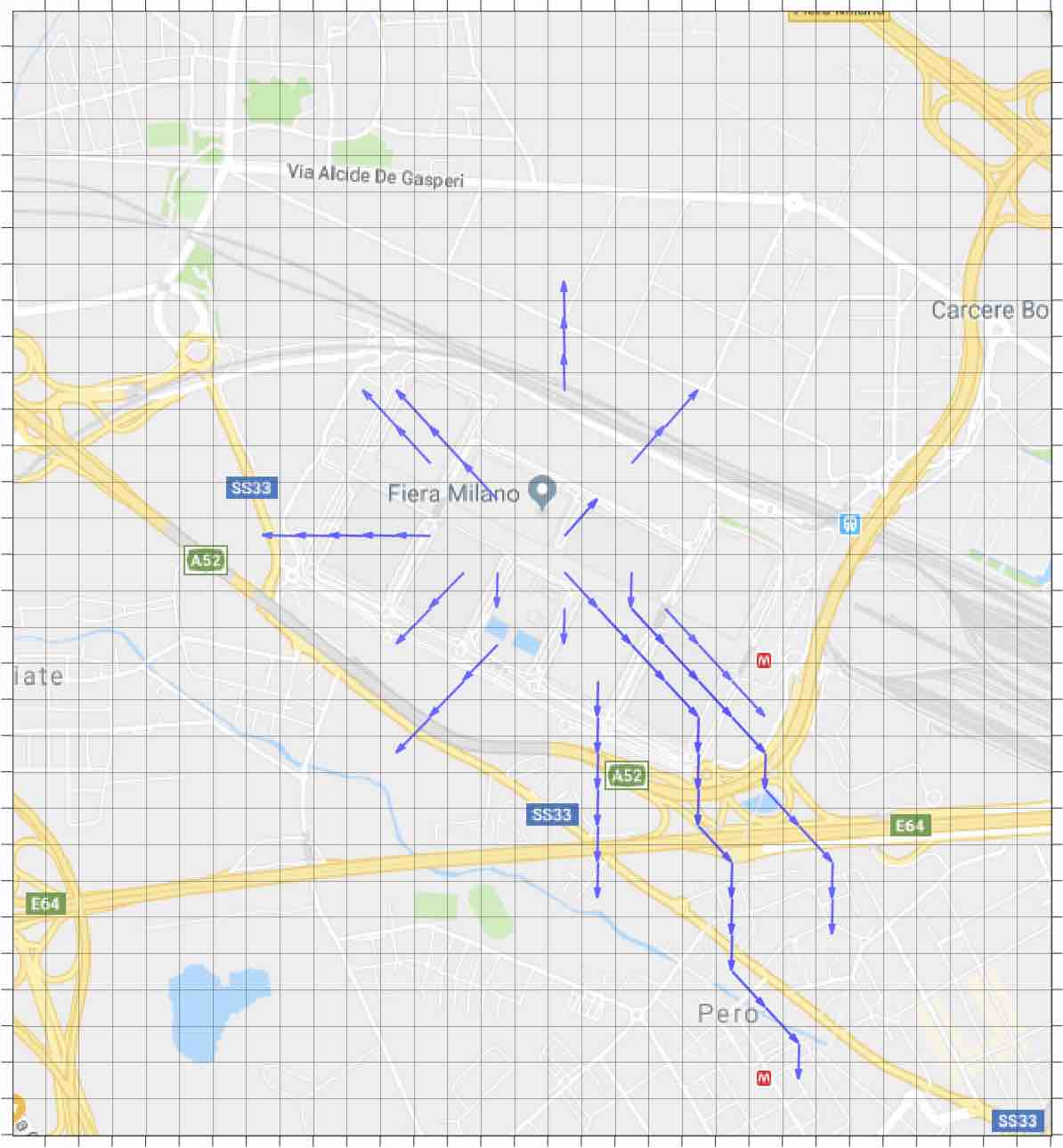}
\caption{Main flows directed to the exhibition area in the morning (09:45 - 09:46) (left), Main flows from the exhibition area in the evening (18:00 - 18:01).(right)}
\label{fig:fiera}
\end{figure}

In this test we show how the algorithm is able to capture the way a big event influences human flows. The event we have considered is the exhibition of the \emph{Salone del Mobile}, held every April at Fiera Milano exhibition center in Rho, near Milan. We analyze a square area of $31\times 31$ ETUs centered around Fiera Milano. In the left panel of Figure \ref{fig:fiera} we show the morning fluxes directed to the exhibition area whereas in the right panel we show the evening flows directed from the exhibition area to the outside. 

\paragraph{CPU times}
For a simulation of 18 hours, from 06:00 to 23:45, on an area of $31\times 31$ ETUs, the local algorithm requires 6 minutes of CPU time while the global approach requires 30 minutes. Furthermore, we observe that the local algorithm works with $1065$ snapshots of data, whereas the global one with $266$ snapshots.

\section{Conclusions} \label{sec:conclusions}
In this paper we have proposed an efficient method to solve an inverse mass transfer problem, consisting in recovering the dynamics underlying the mass displacement. The proposed algorithm prevents to handle large displacements of the mass, thus saving CPU time and memory allocation with respect to other recently proposed solutions.

The application of the methodology to a real dataset describing the movement of people was also investigated. It is useful to note here that the applicability of the Wasserstein-based methodology was not at all obvious, since it is based on assumptions which are not totally fullfilled. Indeed, the choice of the cost function $c$ does not take into account the fact that people move mainly along roads, are stopped by obstacles, buildings, etc., and in general are not free to move in all directions. Moreover, and most important, the computation of the Wasserstein distance stems from a global optimization problem in which the mass is considered as a whole. In other words, the optimal transport map $T^*$ is found by minimizing the cost of the displacement of the whole crowd, without any distinction among people, and with no regards about individual targets. Despite this, the results we have obtained (see especially those in Figure \ref{fig:provincia} and \ref{fig:fiera}) are exactly as one can expect, meaning that the method is, overall, robust enough to work well even if the constitutive assumptions are not totally fullfilled.

\appendix 
\section{The DMD for the viscous Burgers' equation}
\label{sec:toyexampleDMD}
Here we propose a numerical experiment using data generated by the 2D viscous Burgers' equation:
\begin{equation}\label{burg}
\begin{cases}
\partial_t y - \varepsilon \Delta y + y(\partial_{x_{1}}y+\partial_{x_{2}}y) =0&\quad {\bf x}\in \Omega, t\in[0,T]\\
y({\bf x},t) = 0 &\quad {\bf x}\in \partial\Omega\\
y({\bf x},0)= \sin(\pi x_1)\sin(\pi x_2) &\quad {\bf x}\in\Omega,
\end{cases}
\end{equation}
with ${\bf x}=(x_1,x_2)$, $ \Omega = [0,1]\times[0,1]$, $T=1$ and $\varepsilon = 0.01$.
The numerical approximation of equation \eqref{burg} has been carried out by a finite difference method on a square grid with $N=1600$ nodes, using an adaptive Runge-Kutta method for the time integration with temporal step size $\Delta t>0$.
We refer to \cite{QV94} for further details on the numerical method. In what follows, since we do not know the analytical solution of equation, we consider as reference solution of \eqref{burg} its numerical approximation with $\Delta t = 0.0125$.

To show how DMD reconstructs the model we consider as dataset $\mathcal{X}$ the numerical approximation of \eqref{burg} for the following temporal step sizes $\Delta t = \{0.1, 0.05, 0.025, 0.0125\}$ with snapshots matrices of dimension $N\times 11, N\times 21, N\times 41, N\times 81,$ respectively. 
Once the snapshots have been computed we apply Algorithm \ref{Alg_DMD} to reconstruct the solution corresponding to a fixed $\delta t = 0.0125.$ Here, we aim at testing the capability of the DMD method to approximate the dataset for the missing information.

\begin{figure}[h]
	\centering
	\includegraphics[scale=0.28]{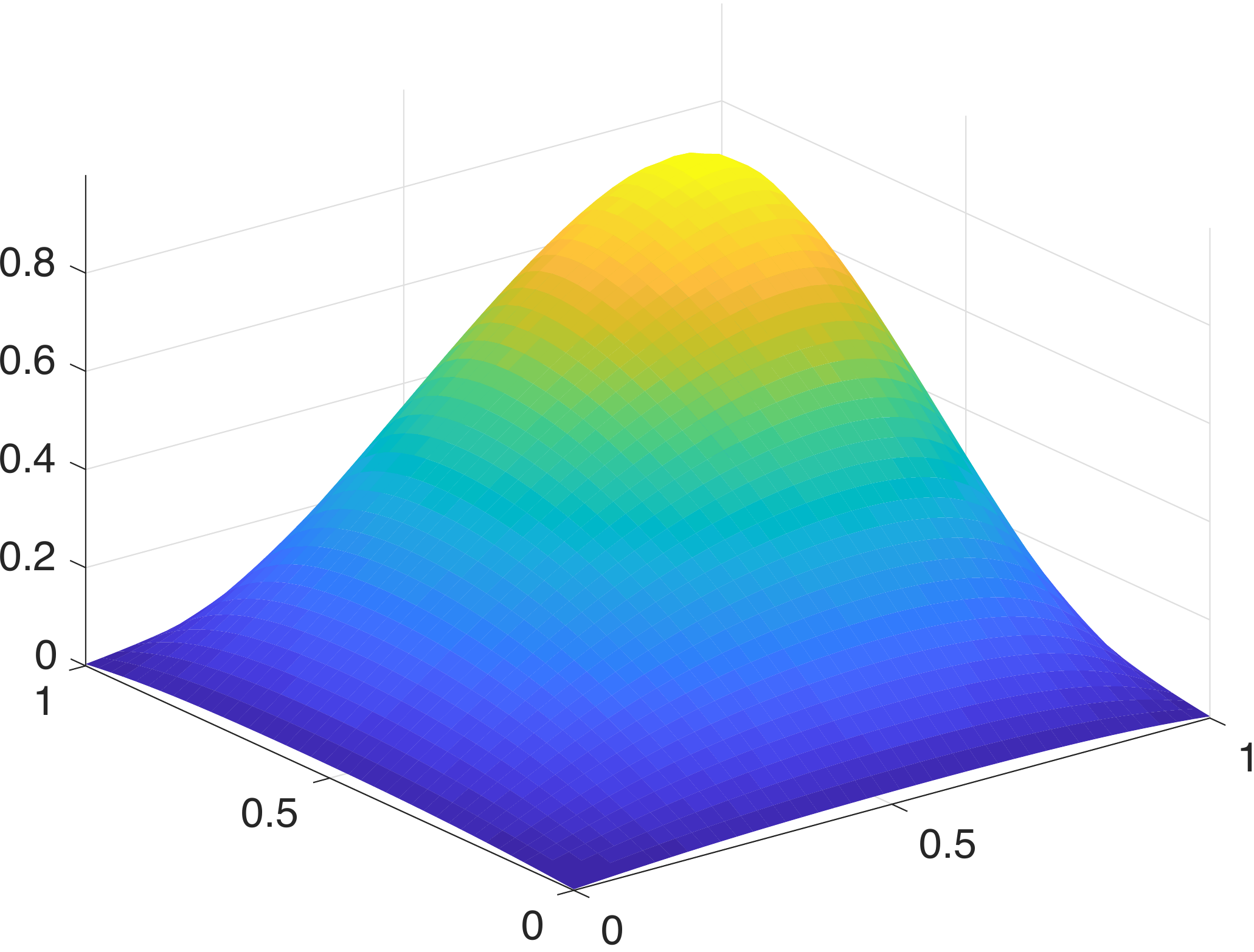}
	\includegraphics[scale=0.28]{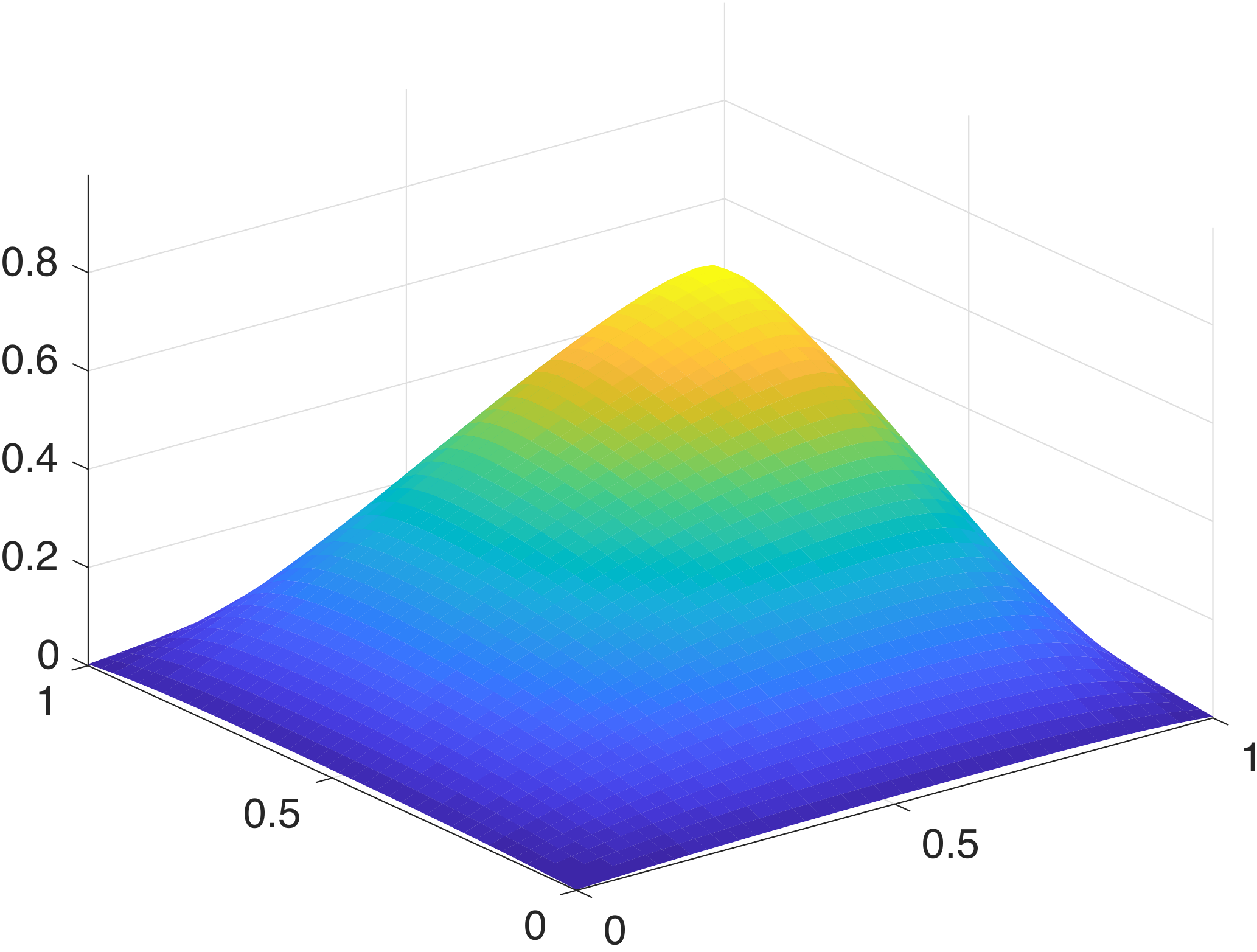}
	\includegraphics[scale=0.28]{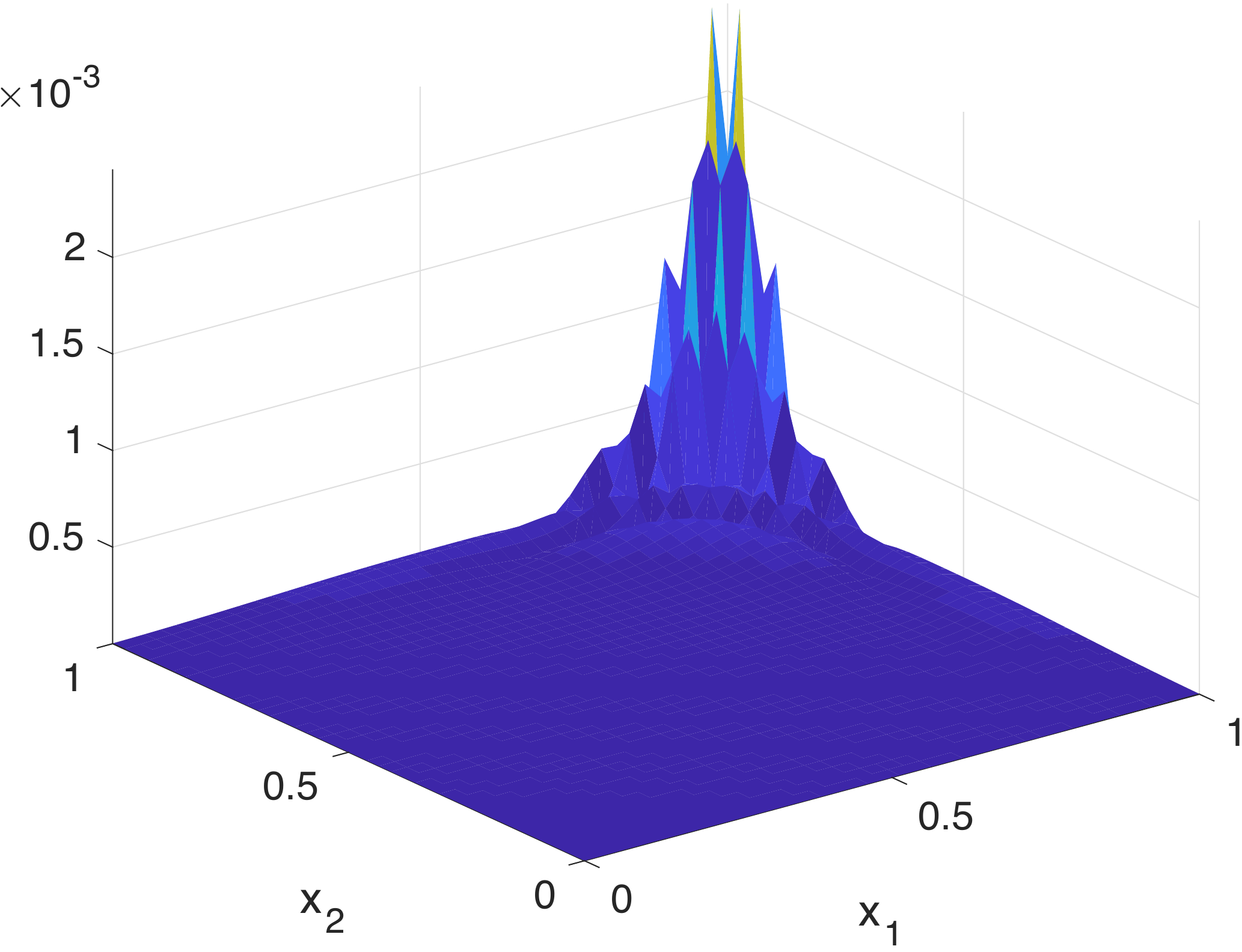}
	\includegraphics[scale=0.28]{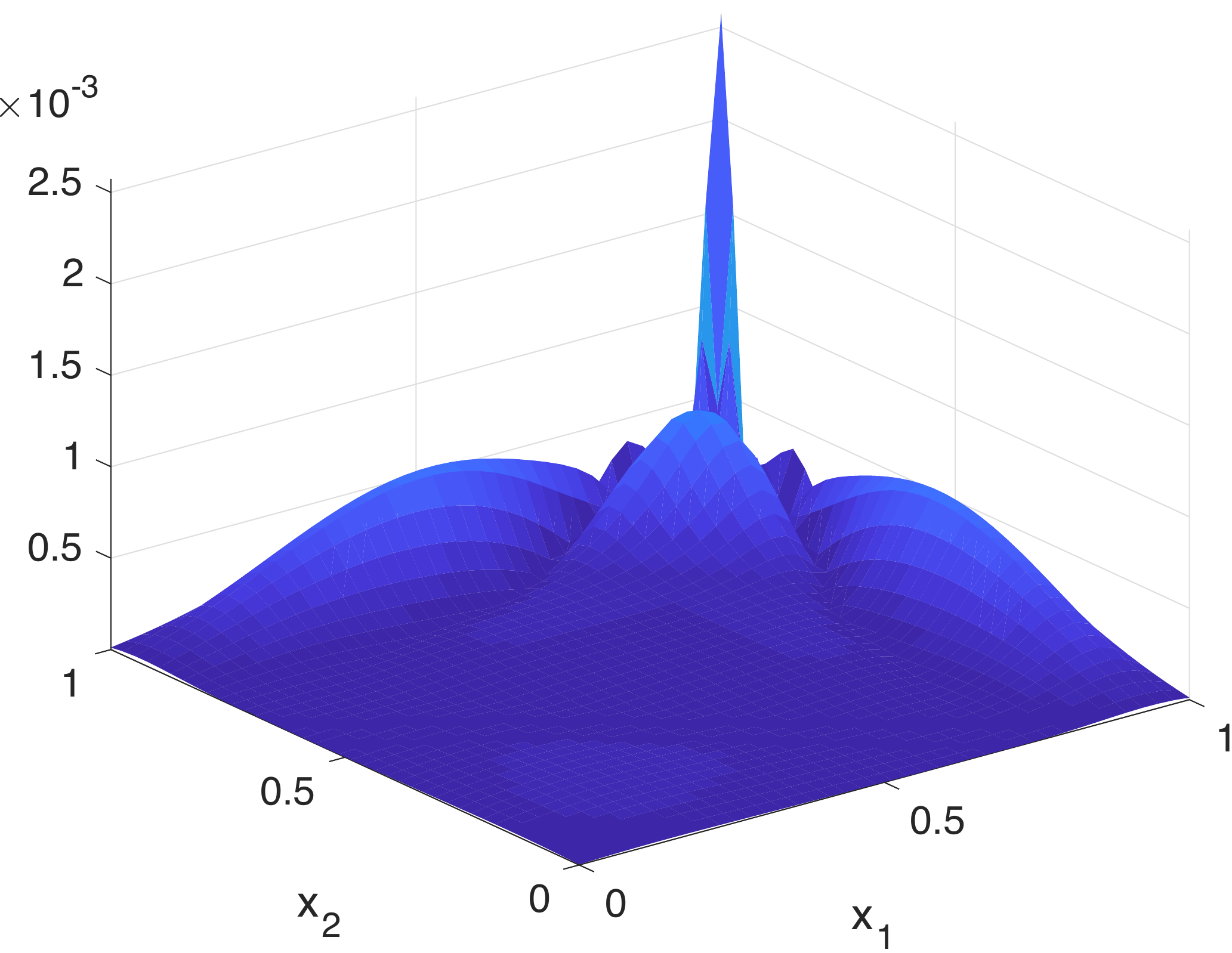}
	\caption{Top: Reference Solution for Burgers' equation at time $t=0.4875$ (left) and $t=0.9875$ (right). Bottom: absolute difference between reference solution and  DMD reconstruction with data $\Delta t =0.1$  at time $t=0.4875$ (left) and $t=0.9875$ (right).}
	\label{fig1:burg}
\end{figure}

In Figure \ref{fig1:burg}, we show the reference solution of \eqref{burg} in the top row and the absolute difference between the reference solution and the reconstructed one with DMD in the bottom row for two time instances. We note that the dataset for the DMD reconstruction corresponds to the solution of \eqref{burg} with $\Delta t = 0.1$. 

\begin{figure}[h!]
	\centering
	\includegraphics[scale=0.23]{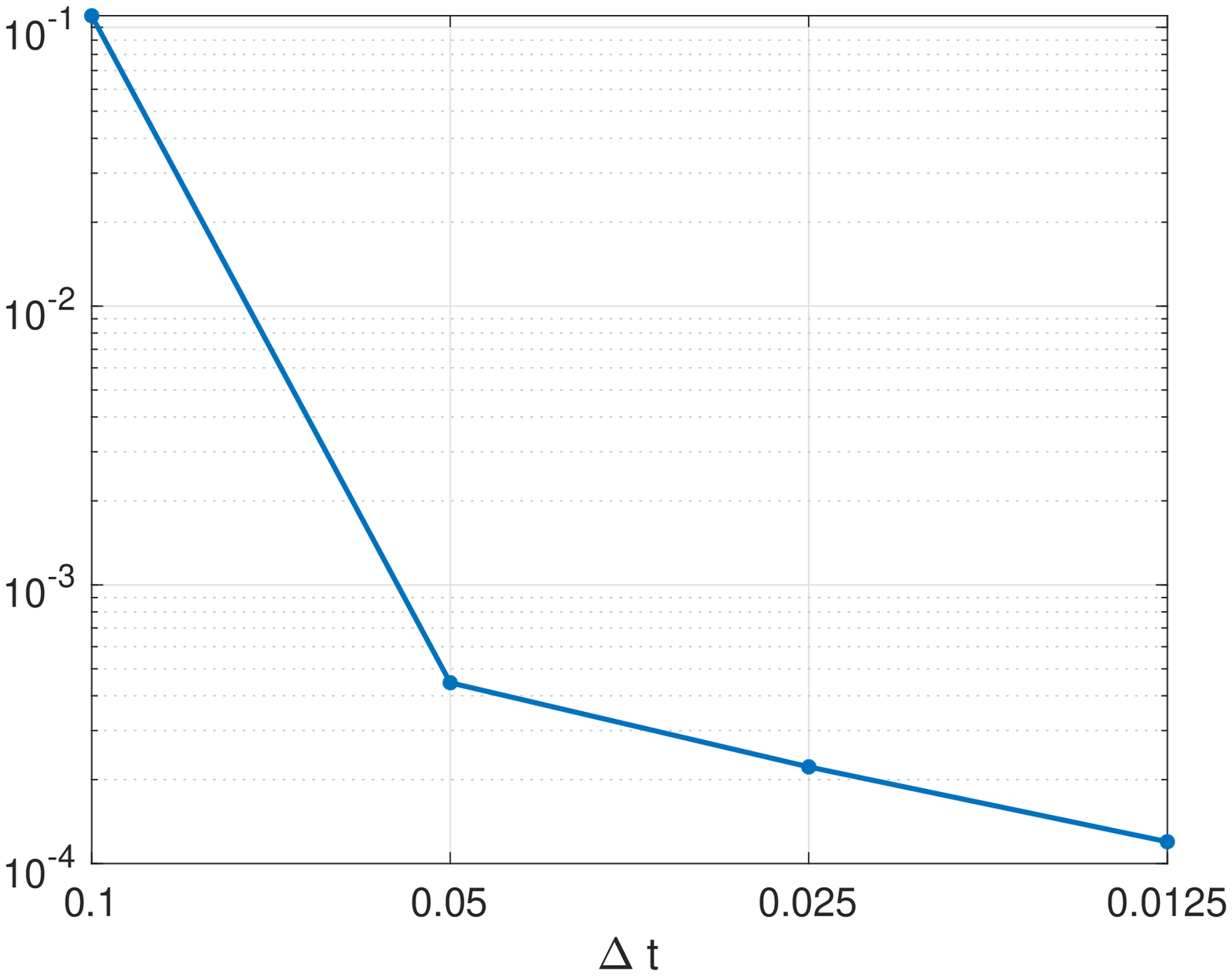}
	\includegraphics[scale=0.23]{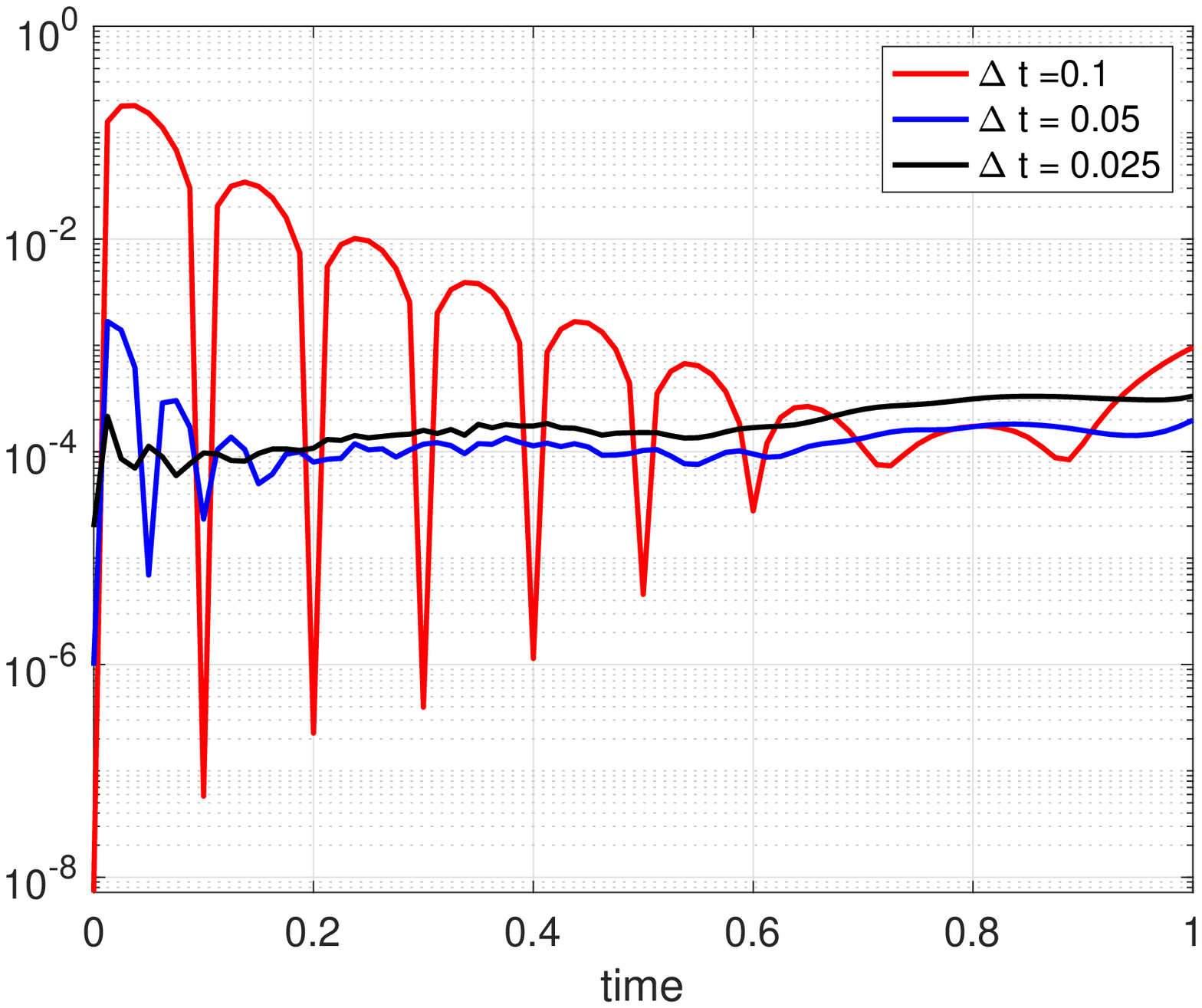}
	\includegraphics[scale=0.29]{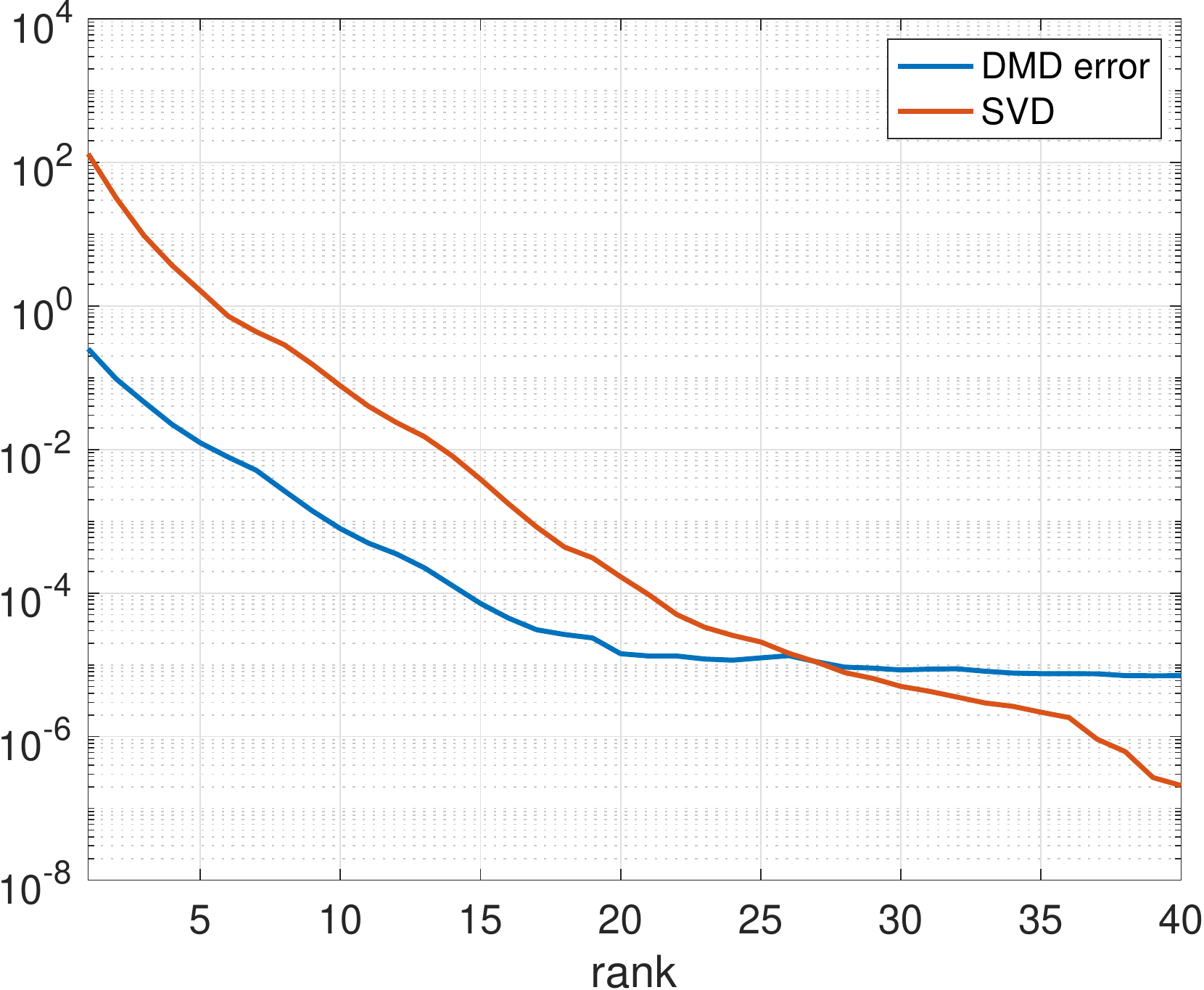}
	\caption{Error analysis of the DMD reconstruction with a reference solution computed with $\Delta t = 0.0125$ and dataset corresponding to $\Delta t = \{0.1, 0.05, 0.025, 0.0125 \}$ as in the $x-$axis using $E^{DMD}$ in \eqref{err:burg} (left), relative error $\widehat{E}^{DMD}(t)$ in \eqref{err:burg} for each time instance for the DMD reconstruction with $\Delta t = \{0.1, 0.05, 0.025\}$ (middle) and singular values of the snapshots matrix with $\Delta t=0.0125$ together with the error behavior $E^{DMD}$ for different rank (right).}
	\label{fig2:burg}
\end{figure}

Then, in Figure \ref{fig2:burg} we show the relative error between the DMD approximation and the reference solution of \eqref{burg} as follows:
\begin{equation}\label{err:burg}
E^{DMD} :=  \frac{\norm{{\bf Y}^R-{\bf Y}^{DMD}}_{F}}{\norm{{\bf Y}^R}_{F}},\qquad
\widehat{E}^{DMD}(t) := \frac{\norm{{\bf y}^R(t)-{\bf y}^{DMD}(t)}_{2}}{\norm{{\bf y}^R(t)}_{2}},\qquad
\end{equation}   
where  ${\bf y}^R(t)$ is the reference solution and ${\bf y}^{DMD}(t)$ is its DMD approximation. The notation ${\bf Y}^R, {\bf Y}^{DMD}$ refers to the matrices where each column contains the reference solution and its DMD approximation, respectively, for different time instances.

 It is clear that the more information on the system, e.g. amount of snapshots, the better error behavior for the DMD solution (see left panel in Figure \ref{fig2:burg}). However, we can already obtain accurate solutions with a poor dataset. In the middle panel of Figure \ref{fig2:burg} we show the error of the model for each time instance for the dataset corresponding to $\Delta t = \{0.1, 0.05, 0.025\}$ using $\widehat{E}^{DMD}(t)$ defined in \eqref{err:burg}. We note that DMD behaves very well for time instances provided in the snapshot set. 
This is to show that DMD is able to reconstruct accurately the solution for missing information. 
In this example the rank $r$ in the DMD method is 9. Finally, we would like to emphasize that, although the DMD model is linear, it works even for nonlinear problem as, e.g., \eqref{burg}. For the sake of completeness we also show in the right panel of Figure \ref{fig2:burg} the decay of error $E^{DMD}$ in \eqref{err:burg} modifying the rank $r$ in the DMD approximation. 
As expected the error decays as the rank increases. The singular values of snapshots set are also in the picture.

\bibliographystyle{plain}
\bibliography{references}

\end{document}